\documentclass{article}

\usepackage{latexsym,amssymb,amsmath}
\usepackage{amsfonts}
\usepackage{color}
\usepackage{centernot}
\usepackage{verbatim}

\usepackage{scalerel}
\usepackage{tikz}
\usetikzlibrary{svg.path}

\definecolor{orcidlogocol}{HTML}{A6CE39}
\tikzset{
  orcidlogo/.pic={
    \fill[orcidlogocol] svg{M256,128c0,70.7-57.3,128-128,128C57.3,256,0,198.7,0,128C0,57.3,57.3,0,128,0C198.7,0,256,57.3,256,128z};
    \fill[white] svg{M86.3,186.2H70.9V79.1h15.4v48.4V186.2z}
                 svg{M108.9,79.1h41.6c39.6,0,57,28.3,57,53.6c0,27.5-21.5,53.6-56.8,53.6h-41.8V79.1z M124.3,172.4h24.5c34.9,0,42.9-26.5,42.9-39.7c0-21.5-13.7-39.7-43.7-39.7h-23.7V172.4z}
                 svg{M88.7,56.8c0,5.5-4.5,10.1-10.1,10.1c-5.6,0-10.1-4.6-10.1-10.1c0-5.6,4.5-10.1,10.1-10.1C84.2,46.7,88.7,51.3,88.7,56.8z};
  }
}

\newcommand\orcidicon[1]{\href{https://orcid.org/#1}{\mbox{\scalerel*{
\begin{tikzpicture}[yscale=-1,transform shape]
\pic{orcidlogo};
\end{tikzpicture}
}{|}}}}

\usepackage{hyperref} 

\def\Gl{\mathop{\rm Gl}\nolimits}
\def\deg{\mathop{\rm deg }\nolimits}
\def\rank{\mathop{\rm rank}\nolimits}

\def\cpr{\mathop{ \; \rm \angle \; }\nolimits}

\newcommand{\se}{\ensuremath{\stackrel{s.e.}{\sim}}}

\newcommand{\FF}{\mathbb F}

\newcommand{\bp}{\mathbf p}
\newcommand{\bq}{\mathbf q}

\newcommand{\bc}{\mathbf c}
\newcommand{\bd}{\mathbf d}

\newcommand{\bu}{\mathbf u}
\newcommand{\bv}{\mathbf v}
\newcommand{\bw}{\mathbf w}

\newcommand{\br}{\mathbf r}
\newcommand{\bs}{\mathbf s}

\newcommand{\ba}{\mathbf a}

\newcommand{\bW}{\mathbf W}

\newcommand{\bz}{\mathbf z}

\newtheorem{theorem}{Theorem}[section]

\newtheorem{lemma}[theorem]{Lemma}
\newtheorem{remark}[theorem]{Remark}
\newtheorem{problem}[theorem]{Problem}
\newtheorem{definition}[theorem]{Definition}
\newtheorem{corollary}[theorem]{Corollary}
\newtheorem{proposition}[theorem]{Proposition}

\providecommand{\abs}[1]{\lvert#1\rvert}

\providecommand{\keywords}[1]
{
  {\small	
  \textbf{Keywords.}} #1
}

\providecommand{\ams}[1]
{
  {\small	
  \textbf{AMS classification codes (2020).}} #1
  }

\title{Bounds for the change of the Weyr characteristic of matrix pencils after 1-rank perturbations
\thanks{This work was 
supported by the Agencia Estatal de Investigaci\'on of Spain MCIN/AEI/10.13039/501100011033/ and by ``ERDF A way of making Europe'' of the ``European Union'' under grants PID2021-124827NB-I00 and RED2022-134176-T.
}
}
\author{Itziar Baraga\~na
\orcidicon{0000-0001-9573-9155}
\thanks{Departamento de Ciencia de la Computaci\'on e I.A., Universidad del Pa\'{\i}s Vasco UPV/EHU, Donostia-San Sebasti\'an, Spain (itziar.baragana@ehu.eus) 
} 
\and 
Alicia Roca
\orcidicon{0000-0000-1926-1895}
\thanks{Corresponding author.}
\thanks{Departamento de Matem\'atica Aplicada, IMM, Universitat Polit\`ecnica de Val\`encia, 46022,Valencia, Spain (aroca@mat.upv.es)
}}

\begin{document}

\maketitle

\begin{abstract}
The complete characterization of the Kronecker structure of a  matrix pencil perturbed by another pencil of rank one is known, and it is stated in terms of very involved conditions.
This paper is devoted to, loosing accuracy,  better understand the meaning of those conditions.

The Kronecker structure of a pencil is determined by the sequences of  
the column and row minimal indices and of the partial multiplicities of the eigenvalues. 
We introduce the Weyr characteristic of a matrix pencil  as the collection of the conjugate partitions of the previous sequences   and provide bounds for the change of each one  of  the Weyr partitions  when the pencil is perturbed by a pencil of rank one.   For each one of the Weyr components, the resulting bounds are expressed only in terms of the corresponding component of the unperturbed pencil.
In order to verify that the bounds are reachable, we also 
characterize the partitions of the  Weyr characteristic of a pencil obtained from another one by removing or adding
 one row. 
 The results hold 
 for  algebraically closed fields.
\end{abstract}


\medskip

\ams{
15A21; 15A22; 15A83; 47A55
}

\section{Introduction}
\label{secintro}

Given a matrix $H$,  the  {\em rank perturbation  problem} consists in characterizing the structure of  $H+P$, where  $P$ is a matrix of  bounded rank. In this work we focus on the 1-rank perturbation problem of matrix pencils, i.e., $H=H(s)$ and $P=P(s)$ are pencils, $\rank(P(s))=1$, and the aim is to  analyze the Kronecker structure  of $H(s)+P(s)$. We recall that the {\em Kronecker structure} of a matrix pencil is formed by the homogeneous invariant factors and the column and row minimal indices, which form a complete system of invariants for the strict equivalence relation. When the pencils are regular, the minimal indices are absent and the Kronecker structure is also called the {\em Weierstrass structure}.

Changes in  the Weierstrass structure of regular   matrix pencils after low rank perturbations have been widely analyzed in the literature in the last decades, both from a generic and  a general point of view. See, for instance 
 \cite{BaRo19, BaRo20b,  BaDoRoSt20,Batzke14, TeDo16, TeDoMo08} and the references therein.
Results for more general pencils have also been obtained: see \cite{DoStSIMAX23} for quasi regular pencils.
For arbitrary pencils the same problem has  been studied  in \cite{TeDo07}, it has been solved 
in \cite{BaRo20,  DoSt20} for $1$-rank perturbations,  and in \cite{DoSt25} for perturbations of the minimal possible rank.

It is worth to stand out the main result in the papers \cite{BaRo20,  DoSt20}, where a solution to the 1-rank perturbation problem  of matrix pencils  is given. In those papers, the complete characterization of the  Kronecker structure of the perturbed pencil  is stated in terms of rather difficult combinatorial conditions. As a consequence, although the conditions can be checked,  it is complicated to know what is the structure of the pencil after perturbation. It seems, therefore, reasonable to try to better  understand the meaning of such characterization, and that is the purpose  of this paper.

A first trial in this direction was done in \cite[Theorem 8.2] {LeMaPhTrWiCAOT21}, where bounds for the changes of certain invariants of a square pencil after some specific type of perturbations of rank one 
(perturbations that does not have row minimal indices) were obtained.
The study was carried out in the
framework of
linear relations. Inspired by this  work, a further step has been done in \cite{BaRo22}. The invariants analyzed in \cite[Theorem 8.2]{LeMaPhTrWiCAOT21} are, in fact, the Weyr characteristic of the regular part of a pencil plus the conjugate partition of its column minimal indices increased in one. In  \cite[Theorem 5.4]{BaRo22} we provided bounds for the changes of these  invariants for arbitrary pencils after arbitrary perturbations of rank one. In the overlapping cases, the achieved bounds improved those of \cite{LeMaPhTrWiCAOT21}. The obtention of the result required to express the main result of \cite{BaRo20} in terms of the conjugate partitions of the sequences defining the Kronecker structure.

At the price of loosing accuracy, the bounds obtained in \cite{BaRo22} display a better understanding of  
the structure of the perturbed pencil than that  derived from \cite{BaRo20}.
But the resulting bounds depend on the complete Weyr characteristics of the original and perturbed  pencils. The aim of  this paper is to 
obtain bounds for the change of the partitions of the Weyr characteristic of the perturbed pencil, expressed only in terms of the corresponding partition of the original pencil.
In order to  verify that the bounds are reachable we solve two 1-row pencil completion problems:
(i) when the Weyr characteristic of a pencil and one of the Weyr invariants of a subpencil is prescribed (Theorems \ref{theo01h} 
  and \ref{theo011h});
(ii) when a subpencil  and one of the invariants of the Weyr characteristic of the completed pencil is prescribed (Theorems \ref{theo01h1} and \ref{theo11h1}). 
As far as we know, the first one of these problems has not been studied yet, while the results in Theorems \ref{theo01h1} and \ref{theo11h1} can be regarded as the conjugate version of some results in \cite{AmBaMaRo24, AmBaMaRo24_2}.

The paper is structured as follows. Section \ref{secpreliminaries} is devoted to introduce notation, algebraic, and combinatorial tools to be used in the paper, and some already known results, either technical or related to pencil completion problems. In Section \ref{bounds} we obtain bounds for the change of the Weyr characteristic of a pencil under  1-row completion and 1-rank perturbation. The two 1-row pencil completion results involving  partial prescription of the Weyr characteristic of the  pencil or of the subpencil mentioned above are  presented in Section \ref{rccompl_ppresc}.  Finally, in Section \ref{reachability} we show with examples that the bounds obtained in Section \ref{bounds} are reachable.

\section{Preliminaries}
\label{secpreliminaries}

All along this paper  $\FF$ is an algebraically closed  field and $\bar \FF=\FF \cup \{\infty\}$. $\FF[s]$ denotes the ring of polynomials in the indeterminate $s$ with coefficients in $\FF$ and  $\FF[s, t]$ is the ring of polynomials in two 
variables $s, t$ with coefficients in $\FF$.
We denote by $\FF^{m\times n}$ and  $\FF[s]^{m\times n}$,  the vector spaces  of $m\times n$ matrices with elements in $\FF$ and  $\FF[s]$  respectively. The degree of a polynomial $p$ is denoted by  $\deg(p)$.
$\Gl_n(\FF)$ is the general linear group of invertible matrices
in $\FF^{n\times n}$. 

A {\em matrix pencil} is a  polynomial matrix $H(s)\in \FF[s]^{m\times n}$ of degree at most one; i.e., $H(s)=H_0+sH_1$, with $H_0, H_1\in \FF^{m\times n}$. The {\em rank} of $H(s)$ is the order of the largest non identically zero minor of $H(s)$, and is denoted by $\rank (H(s))$. It is also known as the {\em normal rank} of $H(s)$.
 If $m=n=\rank(H(s))$, the pencil is called {\em regular}; otherwise it is   {\em singular}.

Two matrix pencils $H(s)=H_0+sH_1, G(s)=G_0+sG_1\in \FF[s]^{m\times n}$ are {\em strictly equivalent}, denoted $H(s)\se G(s)$, 
if there exist invertible matrices $P\in \Gl_m(\FF)$,   $Q\in \Gl_n(\FF)$ such that  $H(s)=PG(s)Q$. When  $H(s)$ and $G(s)$ are not strictly equivalent we write $H(s)\centernot  \se G(s)$.

 A complete system of invariants for  the strict equivalence of a pencil $H(s)\in \FF[s]^{m\times n}$  of $\rank H(s)=\rho$  is formed by the {\em homogeneous invariant  factors} $\pi_1(s,t) \mid \dots \mid \pi_{\rho}(s,t),\ \pi_i(s,t) \in \FF[s,t], \ 1\leq i\leq \rho$,  and two sequences of nonnegative integers $c_1\geq \dots \geq c_{n-\rho}$ and $u_1\geq \dots \geq u_{m-\rho}$, called  the {\em column and row minimal indices} of the pencil, respectively (see  \cite[Ch. 12]{Ga74}).
 They satisfy 
$$\sum_{i=1}^{\rho} d(\pi_i)+\sum_{i=1}^{n-\rho}c_i+\sum_{i=1}^{m-\rho}u_i=\rho.$$
A canonical form for the strict equivalence of matrix pencils is the
{\em Kronecker canonical form}
 (for details see \cite[Ch. 12]{Ga74}).

The transposed pencil $H(s)^T$ has the same homogeneous invariant factors as $H(s)$ and interchanged minimal indices, i.e., the column (row) minimal indices of $H(s)^T$ are the row (column) minimal indices of $H(s)$.
If $H(s)\in \FF[s]^{m\times n}$ and $\rank (H(s))=m$ ($\rank (H(s))=n$), then $H(s)$ does  not have row (column) minimal indices. As a consequence, 
the invariants for the strict equivalence of  regular matrix pencils are reduced to  the homogeneous invariant factors.

The {\em spectrum} of
$H(s)=H_0+sH_1\in \FF[s]^{m \times n}$ is defined as
$
\Lambda(H(s))=\{\lambda\in \overline{\FF}: \rank (H(\lambda))< \rank (H(s))\},
$
where we agree that $H(\infty)=H_1$. The elements $\lambda\in \Lambda(H(s))$ are the {\em eigenvalues} of $H(s)$.

Factorizing  the homogeneous invariant factors of $H(s)$ we can write
$$
\pi_{\rho-i+1}(s,t)=t^{z_i(\infty, H(s))}\prod_{\lambda\in \Lambda(H(s))\setminus\{\infty\}}(s-\lambda t)^{z_i(\lambda, H(s))}, \quad 1\leq i \leq \rho.
$$

The polynomials $\pi_{i}(s,1)$,  $1 \leq i\leq \rho$, are the {\em invariant factors} of $H(s)$, and   
the integers
$
z_1(\lambda, H(s))\geq \dots \geq z_\rho(\lambda, H(s))
$
are called the {\em partial multiplicities} of $\lambda$ in $H(s)$. 
When $z_i(\infty, H(s))>0$, the factor $t^{z_i(\infty, H(s))}$ is an {\em infinite elementary divisor} of the pencil.
For $\lambda \in \overline{\FF}\setminus \Lambda(H(s))$ we take $z_1(\lambda, H(s))=\cdots=z_{\rho}(\lambda, H(s))=0$. We agree that $z_i(\lambda, H(s))=+\infty$ for $i<1$ and $z_i(\lambda, H(s))=0$ for $i>\rho$, for $\lambda \in \overline{\FF}$. We also agree that $\pi_i(s,t)=1$ for $i<1$ and  $\pi_i(s,t)=0$ for $i>\rho$.

\medskip

In what follows we work with sequences of integers. 
We call  {\em partition} of a positive integer $n$ to  an  infinite 
sequence of nonnegative integers $\br = (r_1, r_2, \ldots)$, almost all zero, such that $r_1\geq r_2 \geq  \ldots$, and $\sum_{i=1}^n r_i=n$. 
We also work with finite sequences of integers
$\bc=(c_1, \dots, c_m)$, whith 
$c_1\geq \dots \geq c_m$. When necessary, we take $c_i=+\infty$ if $i<1$ y $c_i=-\infty$ if $i>m$. 

Given a partition $\bc=(c_1, c_2, \dots)$ or a finite sequence  $\bc=(c_1, \dots, c_m)$ of nonnegative integers, we denote by $\abs{\bc}$ the sum of its components, $\abs{\bc}=\sum_{i\geq 1}c_i$, and
the {\em conjugate partition} of  $\bc$ is the partition $\br= (r_1, r_2, \ldots)$ defined
as $r_k := \#\{i : c_i \geq k\}, \ k\geq 1$.
In general, we refer to the conjugate partition of $\bc$ as $\bar \bc$.
For a finite sequence  of nonnegative integers $\bc=(c_1, \dots, c_m)$ we also define  the term $r_0= \#\{i : c_i \geq 0\}$ and  write
$\br_*=(r_0, r_1, r_2, \ldots)$. Notice that $r_0= m$. 
On the other hand, whenever we have a partition $\br_*=(r_0, r_1, r_2, \ldots)$, 
we will take $\br=(r_1, \dots)$.

Given two partitions (finite sequences of integers) $\bc$ and $\bd$, 
$\bc+\bd$ is the partition (finite sequence of integers) whose components are the sums of the components of $\bc$ and $\bd$.

If $a$ and $b$ are integers such that $a>b$ we agree that $\sum_{j=a}^{b}=0$.
Given a real number $x$, we denote by $[x]$  the integer part of $x$.

Let $H(s)\in \FF[s]^{m\times n}$ be a pencil of $\rank H(s)=\rho$ having $\bz(\lambda)=(z_1(\lambda),z_2(\lambda), \dots )$  as the partition of 
partial multiplicities at $\lambda\in \bar \FF$     and $\bc=(c_1, \dots, c_{n-\rho})$ and $\bu=(u_1, \dots, u_{m-\rho})$ as the sequences of column and row minimal indices, respectively{. Recall that we take   $z_1(\lambda)= z_2(\lambda)=\dots=0$ if $\lambda \not \in \Lambda(P(s))$.
The {\em Weyr characteristic  at $\lambda\in \bar \FF$ of $H(s)$}   is  the conjugate partition of $\bz(\lambda)$, and is denoted by $\bw(\lambda)$, i.e., $\bw(\lambda)=\bar \bz(\lambda)=(w_1(\lambda), w_2(\lambda), \dots )$. Moreover, if
 $\Lambda(P(s))=\{\lambda_1, \dots, \lambda_\ell\}$ with $\ell\geq 0$, we define  the {\em Weyr characteristic of the regular part of $H(s)$} as the collection of partitions $\bW=(\bw(\lambda_1), \dots, \bw(\lambda_\ell))$. Observe that if $\ell=0$, then $\bW=((0),\dots)$.
The {\em column and row Brunovsky partitions} of $H(s)$ are  two partitions $\br_*=(r_0, r_1, \dots)$ and $\bs_*=(s_0, s_1, \dots)$ respectively, where $\br=\bar \bc=(r_1,r_2, \ldots)$ is the conjugate partition of $\bc$, $\bs =\bar \bu=(s_1,s_2, \ldots)$ is the 
 conjugate partition of $\bu$, $r_0=n-\rho$ is the number of column minimal indices and $s_0=m-\rho$
is the number of row minimal indices. 
The {\em Weyr characteristic of $H(s)$} is formed by the collection of partitions $(\bW, \br_*,  \bs_*)$, where $\bW$, $\br_*$ and $\bs_*$ are, respectively, the  Weyr characteristic of the regular part and  the column and row Brunovsky partitions of $H(s)$. We would like to emphasize that the Weyr characteristic of a pencil  is a complete system of invariants for the strict equivalence relation.

\medskip

The main target of this work is to solve the next problem. It was proposed by Prof. Carsten Trunk in terms of linear relations. Here we pose it and give a solution for matrix pencils.

\begin{problem}\label{probpert}
    Let  $H(s)\in \FF[s]^{m\times n}$ be a  matrix pencil with Weyr characteristic $(\bW, \br_*, \bs_*)$ and let $P(s)\in \FF[s]^{m\times n}$ be a matrix pencil of
    $\rank (P(s))=1$. If
      $(\widehat \bW, \widehat \br_*, \widehat \bs_*)$   is the Weyr characteristic of the perturbed pencil $H(s)+P(s)$,  find bounds for the components of $\widehat \bW$ ($\widehat\br_*$, $\widehat \bs_*$) in terms of those of  $\bW$ ($\br_*$, $\bs_*$).  
  \end{problem}
In other words, if we know the Weyr characteristic of the regular part  of a given pencil, and we perturb it with a pencil of rank one, the problem consists in finding bounds for the possible Weyr characteristics of the regular part  of the perturbed pencil. And the same considerations for the column and row partitions.

  \medskip

In the following, we present a number of results about perturbation and completion of pencils that will be used later.

The next lemma is a rephrasement of \cite[Proposition 2.1]{BaRo20} (see also \cite[Proposition 3.1]{GeTr17}).
  
  \begin{lemma}{\rm \cite[Proposition 2.1]{BaRo20}}
    \label{pencilrank1}
      Let $P(s)\in \FF[s]^{m\times n}$ be a matrix pencil of
      $\rank (P(s))=1$.
      \begin{enumerate}
\item \label{ituvs}
If $P(s)$ does  not have positive row minimal indices, then there exist a nonzero vector $u\in \FF^m$ and a nonzero  (possibly constant) pencil $v(s)\in \FF[s]^n$ such that  $P(s)=uv(s)^T$.
\item \label{itusv}
If $P(s)$ does  not have  positive column minimal indices, then there exist a nonzero  (possibly constant) pencil $u(s)\in \FF[s]^m$ and a nonzero vector $u\in \FF^n$ such that  $P(s)=u(s)v^T$.
        \end{enumerate}
        
\end{lemma}

The 1-rank perturbation problem of a pencil can be related to the pencil completion problem as follows.

\begin{lemma}{\rm \cite[Lemma 3.3 and Remark 3.4]{BaRo20}}
\label{lemmaeq}
Let  $H(s), G(s)\in \FF[s]^{m \times n}$ be matrix pencils such that 
$H(s)\not \se G(s)$.
\begin{enumerate}
\item \label{itpertrow}
  There exists a matrix pencil  $P(s)=uv(s)^T \in\FF[s]^{m\times n}$ such that
$H(s)+P(s)\se G(s)$  if and only if there exist matrix pencils  $h(s), g(s)\in\FF[s]^{1\times n}$, $H_{1}(s)\in \FF[s]^{(m-1)\times n}$ such that
$ H(s)\se\begin{bmatrix}h(s)\\H_1(s)\end{bmatrix}$ and
$G(s)\se\begin{bmatrix}g(s)\\H_{1}(s)\end{bmatrix}$.
\item \label{itpertcol}
  There exists a matrix pencil  $P(s)=u(s)v^T \in\FF[s]^{m\times n}$ such that
$H(s)+P(s)\se G(s)$  if and only if  there exist matrix pencils $h'(s), g'(s)\in\FF[s]^{m\times 1}$, $H'_{1}(s)\in \FF[s]^{m\times (n-1)}$ such that
$ H(s)\se\begin{bmatrix}h'(s)&H'_1(s)\end{bmatrix}$ and
$G(s)\se\begin{bmatrix}
g(s)&H'_1(s)\end{bmatrix}$.
\end{enumerate}
\end{lemma}

The  row completion pencils  results below (Lemmas \ref{lemmaDox0} and \ref{lemmaDox1}) involve a combinatorial object known as 1step-majorization that is a particular case of the generalized majorization introduced in \cite[Definition 2]{DoStEJC10}.

\begin{definition}
Given two finite sequences of  integers $\bc = (c_1, \dots, c_p)$ and  $\bd=(d_1, \dots, d_{p+1})$, 
it is said that  $\bd$ is {\em 1step-majorized} by $\bc$   (denoted $\bd \prec' \bc$) if
$$c_{i}=d_{i+1}, \quad h \leq i \leq p, $$
where $h=\min\{i: c_i<d_i\}$ ($c_{p+1}=-\infty$).
\end{definition}

Lemmas  \ref{lemmaDox0} and \ref{lemmaDox1} are particular cases of \cite[Theorem 4.3]{DoSt19} and they can also be found in  \cite[Lemmas 4.3 and 4.4]{BaRo20} (see also   \cite[Lemma 5.10]{BaRo20}).

\begin{lemma}
\label{lemmaDox0}
Let $H(s)\in \FF[s]^{ m\times n}$, $H_1(s)\in \FF[s]^{( m-1)\times n}$ be two matrix pencils,  $\rank(H_1(s))=\rank(H(s))=\rho$. Let 
$\pi^1_1(s, t)\mid \dots \mid \pi^1_{\rho}(s, t)$,
$\bc=(c_1, \dots,c_{n- \rho})$ and
$\bu=(u_1, \dots, u_{m-\rho-1})$
be the homogeneous invariant factors  and  the sequences of column and row minimal indices of $H_1(s)$, respectively, and let 
$\pi_1(s, t)\mid \dots \mid \pi_{\rho}(s, t)$, 
$\bd=(d_1,  \dots, d_{n-\rho})$ and
$\bv=(v_1, \dots, v_{m-\rho})$
be the homogeneous invariant factors and  the  sequences of column and  row minimal indices  indices of $H(s)$, respectively.

There exists a pencil $h(s)\in\FF[s]^{1 \times n}$ such that $H(s)\se \begin{bmatrix}
  h(s)\\H_1(s)\end{bmatrix}$ if and only if
\begin{equation}\label{inter}
\pi_i(s, t)\mid\pi^1_i(s, t)\mid\pi_{i+1}(s, t), \quad 1\leq i \leq \rho, 
\end{equation}
\begin{equation*}\label{coleq}
\bd= \bc,
\end{equation*}
\begin{equation*}\label{rowprec}
\bv\prec'\bu.
\end{equation*}
\end{lemma}

\begin{lemma}
  \label{lemmaDox1}
  Let $H(s)\in \FF[s]^{ m\times n}$, $H_1(s)\in \FF[s]^{( m-1)\times n}$ be two matrix pencils,  $\rank(H(s))=\rho$, $\rank(H_1(s))=\rho-1$.
  Let 
$\pi^1_1(s, t)\mid \dots \mid \pi^1_{\rho-1}(s, t)$,
$\bc=(c_1, \dots, c_{n-\rho+1})$ and
$\bu=(u_1, \dots, u_{m-\rho})$
be the homogeneous invariant factors and the sequences of column and  row minimal indices of $H_1(s)$, respectively, and let 
$\pi_1(s, t)\mid \dots \mid \pi_{\rho}(s, t)$, 
$\bd=(d_1,\dots, d_{n-\rho})$ and
$\bv=(v_1, \dots,  v_{m-\rho})$
be the homogeneous invariant factors and the sequences of column and  row minimal indices  indices of $H(s)$, respectively.

There exists a pencil $h(s)\in\FF[s]^{1 \times n}$ such that $H(s)\se \begin{bmatrix}
  h(s)\\H_1(s)\end{bmatrix}$ if and only if
(\ref{inter}), 
\begin{equation*}\label{colprec}
\bc\prec' \bd,
\end{equation*}
\begin{equation*}\label{rowequal}
\bv=\bu.
\end{equation*}

\end{lemma}

We have to express these results in terms of the Weyr characteristics of the pencils involved. It requires to define ``the conjugate version'' of the 1step-majorization. This majorization, without assigning it a specific name, was first introduced in \cite{BaZa90}. 
 
 \begin{definition}\label{defmgc}
Given two partitions $\br_* = (r_0, r_1, \dots)$ and  
$\bs_*=(s_0, s_1, \dots)$ 
it is said that  $\bs_*$ is {\em conjugate majorized}  by $\br_*$ (denoted $\bs_* \cpr \br_*$) 
if $r_0=s_0+1$ and 
$$r_{i}=s_i+1, \quad 0 \leq i \leq g, $$
where $g=\max\{i: r_i>s_i\}$.
     \end{definition}

The subsequent Propositions \ref{propDox0conj} and \ref{propDox1conj} establish Lemmas \ref{lemmaDox0} and \ref{lemmaDox1} in terms of the Weyr characteristics.  
The proof is an immediate  consequence of the following auxiliary lemmas.
\begin{lemma}{\rm \cite[Lemma 3.2]{LiSt09}, \cite[Lemma 4.3]{BaRo22} }
\label{lemmaconj}
Let $(a_1, a_2, \dots, )$ and $(b_1, b_2, \dots, )$ be two partitions. Let $(p_1, p_2, \dots)=\overline{(a_1,a_2, \dots, )}$ and $(q_1, q_2, \dots)=\overline{(b_1, b_2, \dots, )}$ be their conjugate partitions.
Let $k\geq 0$ be an integer. 
Then, $a_{j}\geq b_{j+k}$ for   $j\geq 1$
if and only if 
$
p_{j}\geq q_{j}-k$ for  $j\geq 1$.
\end{lemma}
\begin{lemma}{\rm\cite[Proposition 4.5]{BaRo22} }
  \label{propconj}
Given two finite sequences of nonnegative  integers $\bc = (c_1, \dots, c_{m+1})$ and  
$\bd=(d_1, \dots, d_{m})$, let
$\br=(r_1, \dots)=\overline{\bc}$,
$\bs=(s_1, \dots)=\overline{\bd}$ be their conjugate partitions, 
$r_0=m+1=s_0+1$, and
$\br_*=(r_0, r_1, \dots)$, $\bs_*=(s_0, s_1, \dots)$. 
Then $\bc\prec'\bd$ if and only if  $\bs_* \cpr \br_*$.
\end{lemma}

\begin{proposition}
  \label{propDox0conj}
  Let $H_1(s)\in \FF[s]^{(m-1)\times n}$, $H(s)
  \in \FF[s]^{ m\times n}$ be two matrix pencils,
  $\rank(H(s))=\rank(H_1(s))$. Let 
$(\bW^1, \br_*^1,  \bs_*^1)$ and $(\bW, \br_*,  \bs_*)$ be the  Weyr characteristics of $H_1(s)$ and  $H(s)$, respectively.
There exists a pencil $h(s)\in\FF[s]^{1 \times n}$ such that $H(s)\se \begin{bmatrix}
  h(s)\\H_1(s)\end{bmatrix}$ if and only if
\begin{equation}\label{interww1w+1}
w_i(\lambda)\leq  w^1_i(\lambda)\leq  w_i(\lambda)+1, \quad i\geq 1, \quad \lambda \in \bar \FF,
\end{equation}
\begin{equation}\label{coleqconj}
\br_*= \br^1_*,
\end{equation}
\begin{equation}\label{rowprecconj}
\bs_*^1\cpr\bs_*.
\end{equation}
\end{proposition}

\begin{proposition}
  \label{propDox1conj}
  Let $H_1(s)\in \FF[s]^{(m-1)\times n}$, $H(s)
  \in \FF[s]^{ m\times n}$ be two matrix pencils,
  $\rank(H(s))=\rank(H_1(s))+1$. Let 
$(\bW^1, \br_*^1,  \bs_*^1)$ and $(\bW, \br_*,  \bs_*)$ be the  Weyr characteristics of $H_1(s)$ and  $H(s)$, respectively.
There exists a pencil $h(s)\in\FF[s]^{1 \times n}$ such that $H(s)\se \begin{bmatrix}
  h(s)\\H_1(s)\end{bmatrix}$ if and only if
\begin{equation}\label{interw-1w1w}
w_i(\lambda)-1\leq  w^1_i(\lambda)\leq  w_i(\lambda), \quad i\geq 1, \quad \lambda \in \bar \FF,
\end{equation}
\begin{equation}\label{colprecconj}
\br_*\cpr\br^1_*,
\end{equation}
\begin{equation}\label{roweqconj}
  \bs_*=\bs_*^1.
  \end{equation}
 
\end{proposition}

Finally, we present some technical results.

\begin{remark}\label{remranks}
  With  the  notation of Propositions \ref{propDox0conj} and \ref{propDox1conj},
   assume that
 $H(s)\se \begin{bmatrix}
  h(s)\\H_1(s)\end{bmatrix}$.  
  \begin{itemize}
    \item\label{itremrankseq}
  If $\rank(H(s))=\rank(H_1(s))$, then
  $\abs{\bW}+\abs{\br}+\abs{\bs}=\abs{\bW^1}+\abs{\br^1}+\abs{\bs^1}$.
  From (\ref{interww1w+1}) and (\ref{coleqconj}) we obtain
  $\abs{\bW}\leq \abs{\bW^1}$ and $\abs{\bW} + \abs{\bs}=\abs{\bW^1}+\abs{\bs^1}$.
  Therefore, $\abs{\bs^1}\leq\abs{ \bs}$. By (\ref{rowprecconj}),
  $\abs{\bs_*}=\abs{ \bs}+s_0=\abs{ \bs}+s^1_0+1\geq \abs{ \bs^1}+s^1_0+1=\abs{ \bs_*^1}+1$; i.e., $\abs{\bs_*}\geq \abs{ \bs_*^1}+1$.

  \item \label{itremranksdif}
    In the same way, if $\rank(H(s))=\rank(H_1(s))+1$, then $\abs{\bW} + \abs{\br}=\abs{\bW^1}+\abs{\br^1}+1$, hence
    $\abs{\br}\leq \abs{\br^1}+1$ and  $\abs{\br_*}\leq \abs{\br_*^1}$.
  \end{itemize}
  \end{remark}

 \begin{lemma}\label{lemmag+1}
  Given two partitions 
  $\br_* = (r_0,r_1,  \dots)$ and  
$\bs_*=(s_0, s_1,\dots)$ 
such that $\bs_* \cpr \br_*$,  let
 $
k=\abs{\br}-\abs{\bs}$

and 
$g=\max\{i: r_i>s_i\}$. Then $k\leq g$ and
$$
0\leq s_i-r_i \leq g-k, \quad i\geq g+1.
$$
  \end{lemma}
{\bf Proof.}
We have
$$
\begin{array}{rll}
  r_i=&s_i+1, & 0\leq i \leq g,\\
r_i\leq &s_i, & g+1\leq i.\\
  \end{array}
$$
Therefore, 
$
k=g+\sum_{j\geq g+1}(r_j-s_j)\leq g
$, and
$$
 s_i-r_i=g-k+\sum_{j= g+1}^{i-1}(r_j-s_j)+ \sum_{j\geq i+1}(r_j-s_j)\leq g-k, \quad i\geq g+1.
$$
 \hfill $\Box$

 \begin{lemma}\label{lemmagc}
  Given two partitions 
  $\br_* = (r_0,r_1,\dots)$ and  
$\bs_*=(s_0, s_1,\dots)$ 
such that $\bs_* \cpr \br_*$,  let $g=\max\{i: r_i>s_i\}$ and let 
   $\bc =(c_1, \dots, )$ be the conjugate partition of  $\br$. Then 
    $$g=c_1 \mbox{ or } g\leq c_2.$$ 
  \end{lemma}
    {\bf Proof.} As  $r_{g}=s_{g}+1\geq 1$,  we obtain that $c_1\geq g$. If $c_1>g$, then
    $1\leq r_{g+1}\leq s_{g+1}\leq s_{g}=r_{g}-1$; therefore $r_g\geq 2$ and $c_2\geq g$.
 \hfill $\Box$
  
\section{Bounds for the change of the Weyr characteristic of a pencil under one row completion and 1-rank perturbation}
\label{bounds}

The target of this section is to study the change of the Weyr characteristic of a matrix pencil $H(s)$ under a perturbation by a matrix pencil $P(s)$  of rank one.  By Lemma \ref{pencilrank1}, $P(s)$ can have two different forms. We first study  the problem under the hypothesis of Lemma \ref{pencilrank1}-\ref{ituvs}, which leads to the one row completion problem stated in Lemma \ref{lemmaeq}-\ref{itpertrow}.

Observe that under the conditions of   Lemma \ref{lemmaeq}-\ref{itpertrow}, we have
$$\max\{\rank(H(s)), \rank(G(s))\}-1\leq \rank(H_1(s))\leq \min\{\rank(H(s)), \rank(G(s))\},$$
i.e., one of the following possibilities holds:
\begin{enumerate}
\item $\rank(G(s))=\rank(  H(s))=\rank (H_1(s))$, 
\item $\rank(G(s))=\rank(  H(s))=\rank (H_1(s))+1$,
  \item  $\rank (G(s))= \rank ( H_1(s))=\rank ( H(s))-1$.
\item  $\rank (H(s))= \rank ( H_1(s))=\rank (G(s))-1$.
   \end{enumerate} 
Similar properties follow if  Lemma \ref{lemmaeq}-\ref{itpertcol} is satisfied. 

\medskip

We start by finding bounds for the change of the Weyr characteristic of a pencil under row completion.

\begin{proposition}
  \label{propgbounds}
  Let $H_1(s)\in \FF[s]^{(m-1)\times n}$, $H(s)
\in \FF[s]^{ m\times n}$, $h(s)\in \FF[s]^{ 1\times n}$ be matrix pencils  such that $H(s)\se \begin{bmatrix}
  h(s)\\H_1(s)\end{bmatrix}$. Let 
$(\bW^1, \br_*^1,  \bs_*^1)$ and $(\bW, \br_*,  \bs_*)$ be the  Weyr characteristics of $H_1(s)$ and  $H(s)$, respectively.
\begin{enumerate}
\item
  If $\rank(H(s))=\rank H_1(s)$, then
 \begin{equation}\label{eqgpboundw01req}
      -1\leq w_i(\lambda)-w^1_i(\lambda)\leq 0, \quad i\geq 1, \quad \lambda \in \bar \FF,
 \end{equation}
 \begin{equation}\label{eqgpboundr01lreq}
      r_i-r_i^1 =0,\quad i\geq 0,
    \end{equation}
    \begin{equation}\label{eqgpbounds01lreq}
     1 - \left [ \sqrt{\abs{\bs_*^1}+1}\right ]\leq s_i-s^1_i\leq 1 , \quad i\geq 0. 
    \end{equation}
  \item
If $\rank(H(s))=\rank H_1(s)+1$, then
 \begin{equation}\label{eqgpboundw01rdif}
      0\leq w_i(\lambda)-w^1_i(\lambda)\leq 1, \quad i\geq 1, \quad \lambda \in \bar \FF,
 \end{equation}
 \begin{equation}\label{eqgpboundr01lrdif}
     -1\leq   r_i-r_i^1 \leq \left [\sqrt{ \abs{\br} } \right ], \quad i\geq 0, 
    \end{equation}
 \begin{equation}\label{eqgpbounds01lrdif}
 s_i-s_i^1 =0\quad i\geq 0. 
    \end{equation}
    
   \end{enumerate} 
\end{proposition}
{\bf Proof.}
\begin{enumerate}
  \item
If $\rank(H(s))=\rank(H_1(s))$, by Proposition \ref{propDox0conj},  (\ref{interww1w+1})-(\ref{rowprecconj}) hold. Conditions
(\ref{eqgpboundw01req})  and (\ref{eqgpboundr01lreq}) are equivalent to   (\ref{interww1w+1}) and (\ref{coleqconj}), respectively.

From (\ref{rowprecconj}), the upper bound in (\ref{eqgpbounds01lreq}) holds,
and by Remark \ref{remranks} we  know that 
    $\abs{\bs^1}\leq\abs{ \bs}$. 
Let $g=\max\{i: s_i>s^1_i\}$.
  By
Lemma \ref{lemmag+1},
  $$
0\leq s^1 _i-s_i\leq g-(\abs{\bs}-\abs{ \bs^1})\leq g, \quad i\geq g+1.
$$
Let $M=\max\{s^1 _i-s_i \; : \; i\geq g+1\}=s^1 _\ell-s_\ell$.
Then $\ell\geq g+1\geq M+1$ and $s^1 _\ell \geq M$.
Thus, 
$$
  \abs{\bs^1_*}=s^1_0+\abs{\bs^1}\geq \sum_{i=0}^{\ell} s^1_i\geq (\ell+1)s^1 _\ell\geq (M+2)M=M^2+2M,  
$$
from where we obtain
$$-1-\left[\sqrt{1+\abs{\bs^1_*}}\right]\leq M\leq -1+\left [\sqrt{1+\abs{\bs^1_*}}\right ].$$
Therefore, $s_i^1-s_i=-1\leq -1+\left [\sqrt{1+\abs{\bs^1_*}}\right ]$ for $0\leq i \leq g$, and $s_i^1-s_i\leq M\leq -1+\left [\sqrt{1+\abs{\bs^1_*}}\right ]$ for $i\geq g+1$; i.e., 
the lower bound in (\ref{eqgpbounds01lreq}) holds.

\item
  If $\rank(H(s))=\rank(H_1(s))+1$, by Proposition \ref{propDox1conj},
  (\ref{interw-1w1w})-(\ref{roweqconj}) hold. Conditions
(\ref{eqgpboundw01rdif})  and (\ref{eqgpbounds01lrdif}) are equivalent to   (\ref{interw-1w1w}) and (\ref{roweqconj}), respectively.

  From (\ref{colprecconj}), the lower bound in (\ref{eqgpboundr01lrdif}) holds, and
   by Remark \ref{remranks} we have 
    $\abs{\br^1}-\abs{ \br}\geq -1$. 
Let $g=\max\{i: r^1_i>r_i\}$.
       By
Lemma \ref{lemmag+1},
  $$
0\leq r_i-r^1_i\leq g-(\abs{\br^1}-\abs{ \br})\leq g+1, \quad i\geq g+1.
$$
Let $M=\max\{r _i-r^1_i \; : \; i\geq g+1\}=r_\ell-r^1_\ell$.
Then $\ell\geq g+1\geq M$ and $r_\ell \geq M$.
Thus, 
$$
  \abs{\br}\geq \sum_{i=1}^{\ell} r_i\geq \ell r_\ell\geq M^2;  
$$
 hence,
$$-\left[\sqrt{\abs{\br}}\right]\leq M\leq \left [\sqrt{\abs{\br}}\right ].$$
Therefore, $r_i-r^1_i=-1\leq \left[\sqrt{\abs{\br}}\right]$ for $0\leq i \leq g$ and $r_i-r^1_i\leq M\leq \left[\sqrt{\abs{\br}}\right]$ for $i\geq g+1$; i.e., 
the upper bound in 
(\ref{eqgpboundr01lrdif}) holds.
\end{enumerate}
\hfill $\Box$

As a consequence, if we do not take into account the rank of the pencils involved, we obtain the following result.
\begin{corollary}
  \label{corgbounds}
  Let $H_1(s)\in \FF[s]^{(m-1)\times n}$, $H(s)
\in \FF[s]^{ m\times n}$, $h(s)\in \FF[s]^{ 1\times n}$ be matrix pencils  such that $H(s)\se \begin{bmatrix}
  h(s)\\H_1(s)\end{bmatrix}$. Let 
$(\bW^1, \br_*^1,  \bs_*^1)$ and $(\bW, \br_*,  \bs_*)$ be the  Weyr characteristics of $H_1(s)$ and  $H(s)$, respectively.
Then
 \begin{equation}\label{eqgpboundw01}
      -1\leq w_i(\lambda)-w^1_i(\lambda)\leq 1, \quad i\geq 1, \quad \lambda \in \bar \FF,
 \end{equation}
 \begin{equation}\label{eqgpboundr01l}
      -1\leq r_i-r^1_i\leq \left[\sqrt{\abs{\br}}\right], \quad i\geq 0,
    \end{equation}
    \begin{equation}\label{eqgpbounds01l}
      1-\left[\sqrt{\abs{\bs_*^1}+1}\right]\leq s_i-s^1_i\leq 1, \quad i\geq 0.
    \end{equation}
\end{corollary}

 Notice that in Lemma \ref{lemmaeq}-\ref{itpertrow} there appears the completion of two pencils.
In the next theorem we obtain bounds for the difference of the Weyr characteristics of two pencils obtained by completion as in  Lemma \ref{lemmaeq}-\ref{itpertrow}. The result is a consequence of Proposition \ref{propgbounds}.

\begin{theorem}
  \label{theogbounds}
  Let  $H(s), G(s)
  \in \FF[s]^{ m\times n}$, $H_1(s)\in \FF[s]^{(m-1)\times n}$, $h(s), g(s)\in \FF[s]^{ 1\times n}$ be matrix pencils  such that  
  $H(s)\se \begin{bmatrix}
  h(s)\\H_1(s)\end{bmatrix}$ and $G(s)\se \begin{bmatrix}
  g(s)\\H_1(s)\end{bmatrix}$. Let 
$(\bW, \br_*,  \bs_*)$ and $(\widehat \bW, \widehat \br_*,  \widehat  \bs_*)$ be the  Weyr characteristics of $H(s)$ and  $G(s)$, respectively.
\begin{enumerate}
\item\label{itrrr1}
  If $\rank(G(s))=\rank(  H(s))=\rank (H_1(s))$, then
 \begin{equation}\label{eqgpboundwrrr1}
      -1\leq \widehat w_i(\lambda)-w_i(\lambda)\leq 1, \quad i\geq 1, \quad \lambda \in \bar \FF,
 \end{equation}
 \begin{equation}\label{eqgpboundrrrr1}
      \widehat  r_i-r_i =0,\quad i\geq 0,
 \end{equation}
    \begin{equation}\label{eqgpboundsrrr1}
      \left. \begin{array}{cll}
        0\leq \widehat  s_i-s_i\leq 1, & i\geq 0,& \mbox{if } \bs=(0),\\
      - \left [ \sqrt{\abs{\bs_*}}\right ]\leq \widehat  s_i-s_i\leq \left [ \sqrt{\abs{\bs_*}}\right ],
      & i\geq 0, & \mbox{if }\bs\neq (0).\\
      \end{array}\right\}
    \end{equation}
    
  \item\label{itrrr1+}
    If $\rank(G(s))=\rank(  H(s))=\rank H_1(s)+1$, then (\ref{eqgpboundwrrr1}),
 \begin{equation}\label{eqgpboundrrrr1+}
     \left. \begin{array}{cll}
        0\leq \widehat  r_i-r_i\leq 2, & i\geq 0,& \mbox{if } \br=(0),\\
      - \left [ \sqrt{\abs{\br}}\right ]-1\leq \widehat  r_i-r_i\leq \left [ \sqrt{\abs{\br}-1}\right ]+2,
      &  i\geq 0, & \mbox{if }\br\neq (0),\\
      \end{array}\right\}
    \end{equation}
    \begin{equation}\label{eqgpboundsrrr1+}
      \widehat  s_i-s_i =0,\quad i\geq 0.
    \end{equation}
  \item \label{itrr+} 
    If $\rank (G(s))= \rank ( H_1(s))=\rank ( H(s))-1$, then 
\begin{equation}\label{eqgpboundwrr+}
      -2\leq \widehat  w_i(\lambda)-w_i(\lambda)\leq 0, \quad i\geq 1, \quad \lambda \in \bar \FF,
 \end{equation}
 \begin{equation}\label{eqgpboundrrr+}
     -\left [ \sqrt{\abs{\br}}\right ] \leq \widehat  r_i-r_i\leq 1,\quad i\geq 0,
 \end{equation}
    \begin{equation}\label{eqgpboundsrr+}
    1- \left [ \sqrt{\abs{\bs_*}+1}\right ]\leq \widehat  s_i-s_i\leq 1, \quad i\geq 0.
    \end{equation}
  \item \label{itr+r}
   If $\rank ( H(s))= \rank ( H_1(s))=\rank ( G(s))-1$, then 
 \begin{equation}\label{eqgpboundwr+r}
      0\leq \widehat  w_i(\lambda)-w_i(\lambda)\leq 2, \quad i\geq 1, \quad \lambda \in \bar \FF,
 \end{equation}
 \begin{equation}\label{eqgpboundrr+r}
     -1\leq \widehat  r_i-r_i\leq \left [ \sqrt{\abs{\br}+1}\right ],\quad i\geq 0,
 \end{equation}
    \begin{equation}\label{eqgpboundsr+r}
   -1\leq \widehat  s_i-s_i\leq  -1+ \left [ \sqrt{\abs{\bs_*}}\right ], \quad i\geq 0.
    \end{equation}
   \end{enumerate} 
\end{theorem}
{\bf Proof.} 
Let $(\bW^1,  \br^1_*,   \bs^1_*)$ be the  Weyr characteristic of $H_1(s)$. 
We can write $\widehat  w_i(\lambda)-w_i(\lambda)=\widehat  w_i(\lambda)-w^1_i(\lambda)+w^1_i(\lambda)-w_i(\lambda)$ for $i\geq 1$ and $\lambda \in \bar \FF$, 
$\widehat  r_i-r_i=\widehat  r_i-r^1_i+r_i^1-r_i$, and $\widehat  s_i-s_i=\widehat  s_i-s^1_i+s_i^1-s_i$ for  $i\geq 0$.

\begin{enumerate}
\item Applying Proposition \ref{propgbounds} we obtain
 $$-1\leq \widehat  w_i(\lambda)-w^1_i(\lambda)\leq 0, \quad 0\leq w^1_i(\lambda)-w_i(\lambda)\leq 1, \quad i\geq 1,\quad \lambda \in \bar \FF, $$ 
$$\widehat  r_i-r_i^1=r_i^1-r_i=0, \quad i\geq 0,$$
  $$1 - \left [ \sqrt{\abs{\bs_*^1}+1}\right ]\leq \widehat  s_i-s^1_i\leq 1, \quad 
    -1\leq s_i^1-s_i\leq \left [ \sqrt{\abs{\bs_*^1}+1}\right]-1, \quad i\geq 0,$$
from where conditions  (\ref{eqgpboundwrrr1}) and  (\ref{eqgpboundrrrr1}) follow straightforward. 
    By Remark \ref{remranks} we have
    $\abs{\bs_*^1}+1\leq \abs{\bs_*}$, therefore  if $\bs\neq (0)$, (\ref{eqgpboundsrrr1}) holds.

    If $\bs=(0)$, from Remark 
\ref{remranks} we know that $\abs{\bs^1}\leq \abs{\bs}=0$, i.e., 
$\bs^1=(0)$. By Proposition
\ref{propDox0conj}, we have 
$\bs^1_* \cpr  \bs_*$ and
$\bs^1_* \cpr \widehat \bs_*$.
Then,  $s_0=\widehat s_0=s_0^1+1$ and $ 0\leq \widehat s_i=\widehat  s_i-s_i=
\widehat s_i-s_i^1\leq 1 $ for $i\geq 1$. 
Therefore, 
(\ref{eqgpboundsrrr1}) also holds.

\item
 Proceeding as in item \ref{itrrr1}, from Proposition \ref{propgbounds} we obtain  (\ref{eqgpboundwrrr1}) and (\ref{eqgpboundsrrr1+}), and  
 \begin{equation} \label{th33_2}
  - \left [ \sqrt{\abs{\br}}\right ]-1\leq \widehat  r_i-r_i\leq \left [ \sqrt{\abs{\widehat  \br}}\right ]+1, \quad i\geq 0.
  \end{equation}
  Moreover,  from Proposition  \ref{propDox1conj}, we have $\br_* \cpr \br^1_*$ and $\widehat  \br_* \cpr \br^1_*$. Then, $r_i^1-r_i\leq 1$ for $i \geq 0$.
Note that $\widehat r_0=r^1_0-1= r_0$.

If $\br=(0)$, then $\widehat  r_i-r_i=\widehat r_i\geq 0$ for $i\geq 1$. Therefore, the lower bound in (\ref{eqgpboundrrrr1+}) holds. Since $\br_* \cpr \br^1_*$,   $\br^1=(\overbrace{1, \dots, 1,}^{g}0) $ for some $g\geq 0$, and as from  Remark \ref{remranks} we have $\abs{\widehat  \br}\leq \abs{ \br^1}+1$, we conclude that  $\widehat  r_1\leq 2$ and $\widehat  r_i\leq 1$ for $i\geq 2$; hence  the upper bound in (\ref{eqgpboundrrrr1+}) also holds. 

Assume now that $\br\neq (0)$. The lower bound of (\ref{eqgpboundrrrr1+})  appears in \ref{th33_2}.
   Let $\widehat  g=\max\{i: r^1_i>\widehat  r_i\}$ and $\widehat  M=\max\{\widehat  r _i-r^1_i: i\geq \widehat  g+1\}=\widehat  r_{\widehat  \ell}-r^1_{\widehat  \ell}$.
Let us see that 
\begin{equation}\label{eqbarMr}
\widehat  M-1\leq  \left [ \sqrt{\abs{\br}-1}\right ].
  \end{equation}
  If $\widehat  M \leq 1$, then (\ref{eqbarMr}) is immediate. Assume that $\widehat  M\geq 2$.
  As in the proof of Proposition \ref{propgbounds}, we get
$\widehat  \ell\geq \widehat  g+1\geq \widehat  M$ and $\widehat  r_{\widehat \ell}\geq \widehat  M$.
Since $\widehat  M-1\leq\widehat  g$ and $\widehat  M-1\leq \widehat  \ell$,  we have
$$r_i^1=\widehat  r_i+1\geq \widehat  r_{\widehat \ell}+1 \geq \widehat  M+1, \quad 1\leq i\leq \widehat  M-1.$$ 
We know that  $r^1_i\leq r_i+1$ for $i\geq 0$, then $r_i\geq \widehat  M$ for $1\leq i \leq \widehat  M-1$. Thus, it results
$$\abs{\br}-1\geq \sum_{i=1}^{\widehat  M-1} r_i-1\geq(\widehat  M-1)\widehat  M-1=(\widehat  M-1)^2+\widehat  M-2\geq (\widehat  M-1)^2,$$
which proves  (\ref{eqbarMr}).
As a consequence,
$$\widehat  r_i-r_i=\widehat  r_i-r^1_i+r^1_i-r_i\leq (\left [\sqrt{\abs{\br}-1}\right ]+1)+1, \quad i\geq 0;$$
i.e., the upper bound in (\ref{eqgpboundrrrr1+}) holds.

  \item By Proposition \ref{propgbounds} we have
 $$-1\leq \widehat  w_i(\lambda)-w^1_i(\lambda)\leq 0, \quad -1\leq w^1_i(\lambda)-w_i(\lambda)\leq 0, \quad i\geq 1,\quad \lambda \in \bar \FF, $$ 
    $$\br_*^1=\widehat  \br_*,  \quad \text{and} \quad  - \left [ \sqrt{\abs{\br}}\right ]\leq r_i^1-r_i\leq 1, \quad i\geq 0,$$
  $$1 - \left [ \sqrt{\abs{\bs_*^1}+1}\right ]\leq \widehat  s_i-s^1_i\leq 1,  \quad i\geq 0, \quad \text{and} \quad \bs_*^1=\bs_*.$$
Then, conditions  (\ref{eqgpboundwrr+}),  (\ref{eqgpboundrrr+}) and (\ref{eqgpboundsrr+}) follow straightforward.

\item
From Proposition \ref{propgbounds} we directly obtain (\ref{eqgpboundwr+r}),
$$
-1\leq \widehat  r_i -r_i^1\leq \left [ \sqrt{\abs{\widehat  \br}}\right ], \quad i\geq 0, \quad \text{and} \quad \br^1_*=\br_*,$$
$$  -1\leq s^1_i-s_i\leq -1 + \left [ \sqrt{\abs{\bs_*^1}+1}\right ], \quad i\geq 0, \quad \text{and} \quad \bs^1_*=\widehat  \bs_*.$$ 
 By Remark  \ref{remranks} we have  $\abs{\widehat  \br}\leq \abs{\br^1}+1=\abs{\br}+1$, and
 $\abs{\bs_*^1}+1\leq\abs{\bs_*}$. Taking into account these inequalities  we obtain  (\ref{eqgpboundrr+r}) and (\ref{eqgpboundsr+r}).
  \end{enumerate}
  \hfill $\Box$

Now we give a solution to Problem \ref{probpert} when the perturbation pencil $P(s)$ does not have  positive row minimal minimal indices; i.e, when $P(s)=uv(s)^T$. The result is a consequence of Lemma \ref{lemmaeq}  and Theorem \ref{theogbounds}.
\begin{theorem}
  \label{theogpertboundsrow}
Let $H(s)\in \FF[s]^{m\times n}$ be a matrix pencil with Weyr characteristic $(\bW, \br_*, \bs_*)$, and  $P(s)=uv(s)^T$ be a matrix pencil with
$u\in \FF^{m\times 1}$ and  $v(s)\in \FF[s]^{n\times 1}$. Let  $(\widehat  \bW, \widehat  \br_*, \widehat  \bs_*)$   be the Weyr characteristic of $H(s)+P(s)$. Then 
 
\begin{enumerate}
\item\label{itpertuvsrrrow}
  If $\rank( H(s)+P(s))=\rank( H(s))$, then (\ref{eqgpboundwrrr1}), (\ref{eqgpboundrrrr1+})  and (\ref{eqgpboundsrrr1}) hold.
  \item \label{itpertuvsrr-row} 
    If $\rank (H(s)+P(s))=\rank ( H(s))-1$, then 
    (\ref{eqgpboundwrr+})--(\ref{eqgpboundsrr+}) hold.
    \item \label{itpertuvsrr+row} 
    If $\rank ( H(s)+P(s))=\rank ( H(s))+1$, then (\ref{eqgpboundwr+r})--
      (\ref{eqgpboundsr+r}) hold.
   \end{enumerate} 
\end{theorem}

By transposition we obtain a  solution
 to Problem \ref{probpert} when the perturbation pencil $P(s)$ does not have positive column  minimal minimal indices; i.e, when $P(s)=u(s)v^T$.

\begin{corollary}
  \label{corgpertboundscol}
Let  $H(s)\in \FF[s]^{m\times n}$ be a matrix pencil with Weyr characteristic $(\bW, \br, \bs)$,  and  $P(s)=u(s)v^T$ be a matrix pencil with
$u(s)\in \FF[s]^{m\times 1}$ and  $v\in \FF^{n\times 1}$. Let  $(\widehat  \bW, \widehat  \br, \widehat  \bs)$   be the Weyr characteristic of $H(s)+P(s)$. Then 
 
\begin{enumerate}
\item\label{itpertuvsrrcol}
  If $\rank( H(s))=\rank( H(s)+P(s))$, then (\ref{eqgpboundwrrr1}),
  \begin{equation}\label{eqgpboundrrrr1col}
\left. \begin{array}{cll}
        0\leq \widehat  r_i-r_i\leq 1, &  i\geq 0,& \mbox{if } \br=(0),\\
      - \left [ \sqrt{\abs{\br_*}}\right ]\leq \widehat  r_i-r_i\leq \left [ \sqrt{\abs{\br_*}}\right ],
      & i\geq 0, & \mbox{if }\br\neq (0).\\
      \end{array}\right\}
    \end{equation}
  \begin{equation}\label{eqgpboundsrrr1+col}
\left. \begin{array}{cll}
        0\leq \widehat  s_i-s_i\leq 2, &  i\geq 0,& \mbox{if } \bs=(0),\\
      - \left [ \sqrt{\abs{\bs}}\right ]-1\leq \widehat  s_i-s_i\leq \left [ \sqrt{\abs{\bs}-1}\right ]+2,
      & i\geq 0, & \mbox{if }\bs\neq (0).\\
      \end{array}\right\}
     \end{equation}
  \item \label{itpertuvsrr-col} 
    If $\rank (H(s)+P(s))=\rank ( H(s))-1$, then
    (\ref{eqgpboundwrr+}),
     \begin{equation}\label{eqgpboundrrr+col}
    1- \left [ \sqrt{\abs{\br_*}+1}\right ]\leq \widehat  r_i-r_i\leq 1, \quad i\geq 0,
    \end{equation}
\begin{equation}\label{eqgpboundsrr+col}
     -\left [ \sqrt{\abs{\bs}}\right ] \leq \widehat  s_i-s_i\leq 1,\quad i\geq 0.
 \end{equation}
   
    \item \label{itpertuvsrr+col} 
      If $\rank ( H(s)+P(s))=\rank ( H(s))+1$, then (\ref{eqgpboundwr+r}),
 \begin{equation}\label{eqgpboundrr+rcol}
   -1\leq \widehat  r_i-r_i\leq  -1+ \left [ \sqrt{\abs{\br_*}}\right ], \quad i\geq 0,
    \end{equation}
 \begin{equation}\label{eqgpboundsr+rcol}
     -1\leq \widehat  s_i-s_i\leq \left [ \sqrt{\abs{\bs}+1}\right ],\quad i\geq 0,
 \end{equation}

   \end{enumerate} 
\end{corollary}
If we do not know which  type of   pencil  $P(s)$ is, we must take into account the possibility that  $P(s)=uv(s)^t$ or $P(s)=u(s)v^T$, and we obtain the following result.
\begin{corollary}
  \label{corgpertboundsrowcol}
Let  $H(s)\in \FF[s]^{m\times n}$ be a matrix pencil with Weyr characteristic $(\bW, \br_*, \bs_*)$, and  let $P(s)$ be a matrix pencil of $\rank (P(s))=1$. Let  $(\widehat  \bW, \widehat  \br_*, \widehat  \bs_*)$   be the Weyr characteristic of $H(s)+P(s)$. Then 
 
\begin{enumerate}
\item\label{itpertuvsrrrowcol}
  If $\rank( H(s)+P(s))=\rank( H(s)$, then (\ref{eqgpboundwrrr1}),
  
\begin{equation}\label{eqgpboundrrrr1gilt}
     \left. \begin{array}{cll}
        0\leq \widehat  r_i-r_i\leq 2, & i\geq 0,& \mbox{if } \br=(0),\\
      \min\{- \left [ \sqrt{\abs{\br}}\right ]-1, -\left [ \sqrt{\abs{\br_*}}\right ]\}\leq \widehat  r_i-r_i\leq \max\{\left [ \sqrt{\abs{\br}-1}\right ]+2,\left [ \sqrt{\abs{\br_*}}\right ]\},
      & i\geq 0, & \mbox{if }\br\neq (0),\\
      \end{array}\right\}
    \end{equation}
\begin{equation}\label{eqgpboundrrrr1gilts}
     \left. \begin{array}{cll}
        0\leq \widehat  s_i-s_i\leq 2, & i\geq 0, & \mbox{if } \bs=(0),\\
      \min\{- \left [ \sqrt{\abs{\bs}}\right ]-1, -\left [ \sqrt{\abs{\bs_*}}\right ]\}\leq \widehat  s_i-s_i\leq \max\{\left [ \sqrt{\abs{\bs}-1}\right ]+2,\left [ \sqrt{\abs{\bs_*}}\right ]\},
      & i\geq 0, & \mbox{if }\bs\neq (0).\\
      \end{array}\right\}
    \end{equation}

  \item \label{itpertuvsrr-colu} 
    If $\rank (H(s)+P(s))=\rank ( H(s))-1$, then
    (\ref{eqgpboundwrr+}),
    \begin{equation}\label{eqgpboundrrr+rowcol}
   \min\{1- \left [ \sqrt{\abs{\br_*}+1} \right ], -\left [ \sqrt{\abs{\br}} \right ]\}\leq \widehat  r_i-r_i\leq 1, \quad i\geq 0,
    \end{equation}
\begin{equation}\label{eqgpboundsrr+rowcol}
   \min\{1- \left [ \sqrt{\abs{\bs_*}+1} \right ], -\left [ \sqrt{\abs{\bs}} \right ]\}\leq \widehat  s_i-s_i\leq 1, \quad i\geq 0.
    \end{equation}
    
    \item \label{itpertuvsrr+colu} 
      If $\rank ( H(s)+P(s))=\rank ( H(s))+1$, then (\ref{eqgpboundwr+r}),
 \begin{equation}\label{eqgpboundrr+rrowcol}
   -1\leq \widehat  r_i-r_i\leq  \max\{-1+ \left [ \sqrt{\abs{\br_*}}\right ],\left [ \sqrt{\abs{\br}+1} \right ]\}, \quad i\geq 0,
    \end{equation}
   \begin{equation}\label{eqgpboundsr+rrowcol}
   -1\leq \widehat  s_i-s_i\leq  \max\{-1+ \left [ \sqrt{\abs{\bs_*}}\right ], \left [ \sqrt{\abs{\bs}+1} \right ]\}, \quad i\geq 0.
    \end{equation}

   \end{enumerate} 
\end{corollary}
If we do not know neither the rank of  $ H(s)+P(s)$ nor the type of   $P(s)$, the result we obtain is the following.

\begin{corollary}
  \label{corgpertboundsrowcol2}
Given  a matrix pencil $H(s)\in \FF[s]^{m\times n}$ with Weyr characteristic $(\bW, \br_*, \bs_*)$ and a matrix pencil  $P(s)\in \FF[s]^{m\times n}$ of $\rank (P(s))=1$,
      let  $(\widehat  \bW, \widehat  \br_*, \widehat  \bs_*)$   be the Weyr characteristic of $H(s)+P(s)$. Then 
  \begin{equation}\label{eqgpboundw}
      -2\leq \widehat  w_i(\lambda)-w_i(\lambda)\leq 2, \quad i\geq 1, \quad \lambda \in \bar \FF,
  \end{equation}

\begin{equation}\label{eqgpboundr}
 \left. \begin{array}{cl}
   \min\{1- \left [ \sqrt{\abs{\br_*}+1} \right ],  -1\}    \leq \widehat  r_i-r_i\leq \max\{2, -1+ \left [ \sqrt{\abs{\br_*}}\right ]\}, \  i\geq 0,
 & \mbox{if } \br=(0),\\
 \min\{- \left [ \sqrt{\abs{\br_*}} \right ],  -1-\left [ \sqrt{\abs{\br}}\right ]\}    \leq \widehat  r_i-r_i\leq \max\{ 2+ \left [ \sqrt{\abs{\br}-1}\right ], 
 \left [ \sqrt{\abs{\br_*}}\right ] 
 \}, \  i\geq 0,
 & \mbox{if } \br\neq(0),\\
     \end{array}\right\}     
   \end{equation}
\begin{equation}\label{eqgpbounds}
 \left. \begin{array}{cl}
   \min\{1- \left [ \sqrt{\abs{\bs_*}+1} \right ],  -1\}    \leq \widehat  s_i-s_i\leq \max\{2, -1+ \left [ \sqrt{\abs{\bs_*}}\right ]\}, \  i\geq 0,
 & \mbox{if } \bs=(0),\\
\min\{ -\left [ \sqrt{\abs{\bs_*}} \right ],  -1-\left [ \sqrt{\abs{\bs}}\right ]\}    \leq \widehat  s_i-s_i\leq \max\{ 2+ \left [ \sqrt{\abs{\bs}-1}\right ], 
 \left [ \sqrt{\abs{\bs_*}}\right ] 
 \}, \  i\geq 0,
 & \mbox{if } \bs\neq(0),\\
     \end{array}\right\}
   \end{equation}

\end{corollary}

\section{Row (column) completion of matrix pencils with some  invariants prescribed}
\label{rccompl_ppresc} 
In order to analyze the reachability of the bounds obtained in Section \ref{bounds}, we present several pencil completion results. On one hand,   we  prescribe the Weyr characteristic of a pencil and one of the Weyr invariants of a subpencil (Theorems  \ref{theo01h}   and \ref{theo011h}). On the other one,  we  prescribe a subpencil and one of the  Weyr invariants of the completed one (Theorems  \ref{theo01h1} and \ref{theo11h1}). 
As mentioned in the Introduction,
to the best of our knowledge,
the first one of these problems has not been studied yet, while the results in Theorems \ref{theo01h1} and \ref{theo11h1} are  the conjugate version of some results in \cite{AmBaMaRo24, AmBaMaRo24_2}.

\medskip
To begin with,  we   analyze  the existence of a pencil with prescribed Weyr characteristic. We introduce first some definitions that will be used throughout the document.
Whenever we are given partitions $\bw(\lambda)=(w_1(\lambda), w_2(\lambda), \dots )$, $\lambda \in \bar \FF$,
we define  $\Lambda(\bw)=\{\lambda \in \bar \FF:  \bw(\lambda) \neq (0)\}$. If $\Lambda(\bw)$ is a finite set, i.e., 
$\Lambda(\bw)=\{\lambda_1, \dots, \lambda_\ell\}$ with $\ell\geq 0$,
we  write $\bW=(\bw(\lambda_1) \dots, \bw(\lambda_\ell))$ and $\abs{\bW}=\sum_{i=1}^\ell \abs{\bw(\lambda_i)}$.
As expected, if $\ell=0$,  then $\bW=((0),\dots)$ and $\abs{\bW} =0$.

\begin{lemma}\label{lemmaexistence}
  Let $\br_*=(r_0, r_1, \dots), \bs_*=(s_0, s_1, \dots)$ and  $\bw(\lambda)$, $\lambda \in \bar \FF$, be partitions.
 If $\Lambda(\bw)$ is a finite set and
  $\abs{\bW}+\abs{\br}+\abs{\bs}=\rho$, then there exits a pencil $H(s)\in \FF[s]^{(\rho +s_0)\times (\rho+r_0)}$ of $\rank(H(s))=\rho$ with Weyr characteristic
 $(\bW, \br_*,  \bs_*)$.
\end{lemma}
{\bf Proof.} Let $\bz(\lambda)=\overline{(\bw(\lambda))}$, $\lambda \in \bar \FF$, $\hat \bc=\overline{\br}=(c_1, c_2, \dots)$ and $\hat \bu=\overline{\bs}=(u_1, u_2, \dots)$.
Then $c_i=0$ for $i>r_1$ and $u_i=0$ for $i>s_1$.
Define
$$\pi_{i}(s, t)=t^{z_{\rho-i+1}(\infty)}\prod_{\lambda\in \Lambda((\bw))\setminus\{\infty\}}(s-\lambda t)^{z_{\rho-i+1}(\lambda)}
, \quad 1\leq i \leq \rho,
$$
$$\bc=(c_1, \dots, c_{r_1}, \dots, c_{r_0}), \quad \bu=(u_1, \dots, u_{s_1}, \dots, u_{s_0}).$$
Let $H(s)$ be a matrix  pencil in  Kronecker canonical form with homogeneous invariant factors $\pi_{1}(s, t), \dots, \pi_{\rho}(s, t)$, column minimal indices
$\bc$, and row minimal indices $\bu$.  Then, $H(s)\in \FF[s]^{(\rho+s_0)\times (\rho+r_0)}$, $\rank(H(s))=\rho$, and the Weyr  characteristic of $H(s)$ is 
$(\bW, \br_*,  \bs_*)$.
\hfill $\Box$

\medskip

The following  technical lemma can be considered as the conjugate version of \cite[Lemma 5.7]{BaRo20}.

\begin{lemma}
  \label{lemmatec1} 
  Given a partition $\bp_*=(p_0, p_1, \dots )$,  let  $\ba =(a_1, \dots )$ be the conjugate partition of $\bp=(p_1, \dots )$ and let $x$ be  an integer.
  There exists a partition $\bq_*=(q_0, q_1, \dots )$ such that $\bq_*\cpr\bp_*$ and $\abs{\bq}=\abs{\bp}-x$
    if and only if  
one of the two following conditions holds.
\begin{itemize}
    \item[\rm (a)]  $p_0=1$ and $x=a_1$.
\item[\rm (b)] $p_0>1$ and
     $$x=a_1 \mbox{ or } x\leq a_2.$$ 
\end{itemize}   
   
  \end{lemma}
    {\bf Proof.} Assume that there exists a partition $\bq_*$ such that $\bq_*\cpr \bp_*$ and $\abs{\bq}=\abs{\bp}-x$.  Then $p_0=q_0+1\geq 1$.
    Let  $g=\max\{i: p_i>q_i\}$. By Lemmas \ref{lemmag+1} and \ref{lemmagc}, we have $x\leq g \leq a_1$.

\begin{itemize}
    \item[(a)] If $p_0=1$, then
    $\bp_*=(\overbrace{1,\dots, 1,}^{a_1+1} 0, \dots )$ and
    $\bq_*=(0, \dots )$, hence
     $x=\abs{\bp}-\abs{\bq}=a_1$.

    \item[(b)] Assume that $p_0>1$.
    Let us see that if $x<a_1$, then $x\leq a_2$.  Assume that $a_2<x<a_1$, then $g\geq x>a_2$, and, by Lemma \ref{lemmagc}, $g=a_1$. Therefore, $p_{g+1}=0$ and
     $g=a_1>x=\sum_{i= 1}^{g}(p_i-q_{i})- \sum_{i\geq g+1}q_{i}= 
    g- \sum_{i\geq g+1}q_{i}$, from where we obtain $q_{g+1}> 0$ and $p_g=q_{g}+1\geq q_{g+1}+1\geq 2$; hence
    $g\leq a_2$, which is a contradiction.
\end{itemize}
Let us  prove the converse.
\begin{itemize}
    \item[(a)] Assume that $p_0=1$ and $x=a_1$. 
    Taking
    $\bq_*=(0, \dots )$, we obtain that 
    $\bq_*\cpr\bp_*$ and $\abs{\bq}=\abs{\bp}-a_1=\abs{\bp}-x$.
\item[(b)] Assume that $p_0>1$, and
$x=a_1$ or $x\leq a_2$.

    If $x=a_1$, define
    $$
    \begin{array}{lll}
      q_i=p_i-1, & 0\leq i \leq a_1,\\
      q_i=p_i=0, & a_1<i.\\
      \end{array}
    $$
If $x\leq a_2$, define
$$
    \begin{array}{lll}
      q_i=p_i-1, & 0\leq i \leq a_2,\\
      q_i=p_i=1, &a_2<i\leq a_1,\\
      q_i=p_i+1=1, &a_1<i\leq a_1+a_2-x,\\
      q_i=p_i=0, & a_1+a_2-x<i.\\     
      \end{array}
    $$
In both cases, $q_0\geq q_1 \geq \dots$, and $\bq_*=(q_0, q_1, \dots )$ satisfies  $\bq_*\cpr \bp*$ and $\abs{\bq}=\abs{\bp}-x$.
\end{itemize}
        \hfill $\Box$

\begin{theorem}
  \label{theo01h} 
  Let $H(s)\in \FF[s]^{ m\times n}$ be a matrix pencil with  Weyr characteristic $(\bW, \br_*,  \bs_*)$, and let  $\bu =(u_1, \dots )$ be  the conjugate partition of $\bs$.
  For $\lambda \in \bar \FF$,
  let  $\bw^1(\lambda)=(w^1_1(\lambda), w^1_2(\lambda), \dots )$ be partitions such that  $\Lambda(\bw^1)$ is a finite set, $\Lambda(\bw^1)=\{\lambda_1, \dots, \lambda_\ell\}$,
  and let $\br_*^1=(r_0^1, r^1_1, \dots)$ and $\bs_*^1=(s_0^1, s^1_1, \dots)$  be partitions. 
 There exist pencils 
 $H_1(s)\in \FF[s]^{(m-1)\times n}$, $h(s)\in \FF[s]^{1\times n}$ such that $H(s)\se \begin{bmatrix}h(s)\\H_1(s)\end{bmatrix}$,
$\rank(H_1(s))=\rank (H(s))$, and   
 \begin{enumerate}
 \item \label{ittheo01hw} 
   $\bW^1=(\bw^1(\lambda_1), \dots, \bw^1(\lambda_\ell))$ is the Weyr characteristic  of the regular part  of $H_1(s)$ if and only if (\ref{interww1w+1}),
    \begin{equation}\label{eqs0geq1}
    s_0\geq 1,
\end{equation}
\begin{equation}\label{eqws1}
\mbox{if $s_0=1$, then }
     \abs{\bW^1}-\abs{\bW}=u_1, 
   \end{equation}
   and
 \begin{equation}\label{eqws}  
 \mbox{if $s_0>1$, then }  
     \abs{\bW^1}-\abs{\bW}=u_1 \mbox{ or } \abs{\bW^1}-\abs{\bW}\leq u_2.
   \end{equation}
 \item \label{ittheo01hr} 
   $\br_*^1$
   is the column Brunovsky partition of $H_1(s)$ if and only if (\ref{coleqconj}) and 
    (\ref{eqs0geq1}).
  
 \item \label{ittheo01hs} 
   $\bs_*^1$
   is the row Brunovsky partition of $H_1(s)$ if and only if (\ref{rowprecconj})
   and
   \begin{equation}\label{abss1leqabss}
     \abs{\bs_*^1}<\abs{ \bs_*}.
\end{equation}
  \end{enumerate}
\end{theorem}
{\bf Proof.}
Assume that there exist pencils 
 $H_1(s)\in \FF[s]^{(m-1)\times n}$, $h(s)\in \FF[s]^{1\times n}$ such that $H(s)\se \begin{bmatrix}h(s)\\H_1(s)\end{bmatrix}$,
$\rank(H_1(s))=\rank (H(s))$, and   $(\bW^1, \br^1_*,  \bs^1_*)$ is the
Weyr characteristic of $H_1(s)$. By Proposition \ref{propDox0conj},
(\ref{interww1w+1})--(\ref{rowprecconj}) hold. Since $\rank(H_1(s))=\rank (H(s))$, we have
$\abs{\bW}+\abs{\br}+\abs{\bs}=\abs{\bW^1}+\abs{\br^1}+\abs{\bs^1}$. From (\ref{interww1w+1}) and (\ref{coleqconj}) we obtain $\abs{\bW}\leq \abs{\bW^1}$ and
$\abs{\bs^1}=\abs{\bs}+\abs{\bW}-\abs{\bW^1}$, respectively. By Lemma \ref{lemmatec1} we obtain 
    (\ref{eqs0geq1})--(\ref{eqws}).
Moreover, $\abs{\bs^1}\leq \abs{\bs}$, and since $s_0=s_0^1+1$ we obtain 
(\ref{abss1leqabss}).

Let us prove the converse.
 \begin{enumerate}
 \item
 Assume that  conditions (\ref{interww1w+1})   and
    (\ref{eqs0geq1})--(\ref{eqws}) hold.
   By Lemma \ref{lemmatec1},
   there exists a partition $\widetilde\bs_*^1=(\widetilde s^1_0, \widetilde s^1_1, \dots )$ such that $\widetilde \bs_*^1 \cpr \bs_*$ and 
   $\abs{\widetilde\bs^1}=\abs{\bs}+\abs{\bW}-\abs{\bW^1}$.
    Define $\widetilde\br_*^1=\br_*$.
   Then $\abs{\bW^1}+\abs{\widetilde \br^1}+\abs{\widetilde\bs^1}=\abs{\bW}+\abs{\br}+\abs{\bs}=\rank(H(s))$. We have
$\rank(H(s))+\widetilde s_0^1=\rank(H(s))+s_0-1=m-1$ and $\rank(H(s))+\widetilde r_0^1=\rank(H(s))+r_0=n$.
   By Lemma \ref{lemmaexistence}
   there exits a pencil $H_1(s)\in \FF[s]^{(m-1)\times n}$ of $\rank(H_1(s))=\rank(H(s))$ with Weyr characteristic $(\bW^1, \widetilde \br^1_*, \widetilde \bs^1_*)$.
 As  condition (\ref{interww1w+1}), $\br_*=\widetilde\br_*^1$ and  $\widetilde \bs_*^1 \cpr \bs_*$ hold, by Proposition \ref{propDox0conj},
there exists a pencil
 $h(s)\in \FF[s]^{1\times n}$ such that $H(s)\se \begin{bmatrix}h(s)\\H_1(s)\end{bmatrix}$.

\item
  Assume that conditions
  (\ref{coleqconj}) and (\ref{eqs0geq1}) hold.
  Define
$$
\begin{array}{ll}
\widetilde s^1_i=s_i-1, &0\leq i \leq u_1,\\
\widetilde s^1_i=s_i=0, &i> u_1,\\
\end{array}
$$
and 
$\widetilde \bs^1_*=(\widetilde s^1_0, \widetilde s^1_1, \dots)$. Then $\widetilde \bs_*^1 \cpr \bs_*$ and  $\abs{\widetilde \bs^1}=\abs{\bs}-u_1$. 

Fix $\lambda_0\in \bar \FF$ and for $\lambda \in \bar \FF$ define a partition $\widetilde \bw^1(\lambda)=(\widetilde w^1_1(\lambda),\widetilde  w^1_2(\lambda), \dots)$ as follows:
$$
\begin{array}{ll}
\widetilde   \bw^1(\lambda) =\bw(\lambda), \quad \lambda\in \bar \FF\setminus\{\lambda_0\},\\
\widetilde   \bw^1(\lambda_0)=\bw(\lambda_0)+\overline{(u_1 )}.
  \end{array}
$$
Recall that $\overline{(u_1)}$ denotes  the conjugate partition of $(u_1)$. 
Then  $w_i(\lambda)\leq \widetilde w^1_i(\lambda)\leq  w_i(\lambda)+1$ for $i\geq 1$ and $\lambda \in \bar \FF$, $\Lambda(\widetilde \bw^1)$ is a finite set and $\abs{\widetilde \bW^1}=\abs{\bW}+u_1$; hence
$\abs{\widetilde \bW^1}+\abs{\br^1}+\abs{\widetilde \bs^1}=\abs{\bW}+\abs{\br}+\abs{\bs}=\rank (H(s))$. 
The result follows as in item \ref{ittheo01hw}.

\item Assume that conditions
  (\ref{rowprecconj}) and
  (\ref{abss1leqabss}) hold.
We have  $s_0-1+\abs{\bs^1}<s_0+\abs{\bs}$, from where $\abs{\bs^1}\leq \abs{\bs}$. Then $x= \abs{\bs}- \abs{\bs^1}\geq 0$.
Define $\widetilde \br_*^1=\br_*$.
 Fix $\lambda_0\in \bar \FF$ and for $\lambda \in \bar \FF$ define a partition $\widetilde \bw^1(\lambda)=(\widetilde w^1_1(\lambda), \widetilde w^1_2(\lambda), \dots)$ as follows:
$$
\begin{array}{ll}
\widetilde   \bw^1(\lambda) =\bw(\lambda), \quad \lambda\in \bar \FF\setminus\{\lambda_0\},\\
 \widetilde  \bw^1(\lambda_0)=\bw(\lambda_0)+\overline{(x )}.
  \end{array}
$$
The result follows as in item \ref{ittheo01hr}. 
   \end{enumerate}
\hfill $\Box$

In  Theorems \ref{theo011h} and \ref{theo01h1}, it plays a role the first term of the conjugate partition of a given one. 
Given a partition $\bw=(w_1, \dots)$, recall that  if $\bz=\bar \bw=(z_1, \dots )$ is the conjugate partition of $\bw$, then  
$z_1=\#\{i\geq 1 : w_i\geq 1\}$.
\begin{theorem}
  \label{theo011h} 
  Let $H(s)\in \FF[s]^{ m\times n}$ be a matrix pencil with  Weyr characteristic $(\bW, \br_*,  \bs_*)$, 
 and for $\lambda \in \bar \FF$ let
$z_1(\lambda)=\#\{i\geq 1 : w_i(\lambda) \geq 1\}$. 
  Let  $\bw^1(\lambda)=(w^1_1(\lambda), w^1_2(\lambda), \dots )$, for $\lambda \in \bar \FF$, be partitions such that  $\Lambda(\bw^1)$ is a finite set, $\Lambda(\bw^1)=\{\lambda_1, \dots, \lambda_\ell\}$,
  and let $\br_*^1=(r_0^1, r^1_1, \dots)$ and $\bs_*^1=(s_0^1, s^1_1, \dots)$  be partitions. 
 There exist pencils 
 $H_1(s)\in \FF[s]^{(m-1)\times n}$, $h(s)\in \FF[s]^{1\times n}$ such that $H(s)\se \begin{bmatrix}h(s)\\H_1(s)\end{bmatrix}$,
$\rank(H_1(s))=\rank (H(s))-1$, and   
 \begin{enumerate}
 \item \label{ittheo11hw} 
   $\bW^1=(\bw^1(\lambda_1), \dots, \bw^1(\lambda_\ell))$ is the Weyr characteristic  of the regular part  of $H_1(s)$ if and only if 
   (\ref{interw-1w1w}) and
   \begin{equation}\label{eqwr+}
     \abs{\bW^1}-\abs{\bW}+1\leq \abs{\br}.
   \end{equation}
 \item \label{ittheo11hr} 
   $\br_*^1$ is the column Brunovsky partition of $H_1(s)$ if and only if (\ref{colprecconj}) and 
   \begin{equation}\label{eqrmins1leqz1}
    0\leq \abs{\br^1}-\abs{ \br}+1\leq \sum_{\lambda\in \bar \FF}z_1(\lambda).
\end{equation}
 \item \label{ittheo11hs} 
   $\bs_*^1$
   is the row Brunovsky partition of $H_1(s)$ if and only if (\ref{roweqconj})
   and
   \begin{equation}\label{eqabswabsrgeq1}
     \abs{\bW}+\abs{ \br_*}>r_0.
\end{equation}
  \end{enumerate}
\end{theorem}
{\bf Proof.}
Assume that there exist pencils 
 $H_1(s)\in \FF[s]^{(m-1)\times n}$, $h(s)\in \FF[s]^{1\times n}$ such that $H(s)\se \begin{bmatrix}h(s)\\H_1(s)\end{bmatrix}$,
$\rank(H_1(s))=\rank (H(s))-1$, and   $(\bW^1, \br^1_*,  \bs^1_*)$ is the
Weyr characteristic of $H_1(s)$. By Proposition \ref{propDox1conj},
(\ref{interw-1w1w})--(\ref{roweqconj}) hold. Since $\rank(H_1(s))=\rank (H(s))-1$, we have
$\abs{\bW}+\abs{\br}+\abs{\bs}-1=\abs{\bW^1}+\abs{\br^1}+\abs{\bs^1}$. From   conditions (\ref{interw-1w1w}) and (\ref{roweqconj}) we obtain $0\geq \abs{\bW^1}-\abs{\bW}=
\abs{\br}-\abs{\br^1}-1$, hence (\ref{eqwr+}) and $0\leq \abs{\br^1}-\abs{\br}+1$ hold. Moreover, since $\abs{\bW}+\abs{\br}=\abs{\bW^1}+\abs{\br^1}+1>0$, (\ref{eqabswabsrgeq1}) follows.

For $\lambda \in \bar \FF$,  let $z_1^1(\lambda)=\#\{i\geq 1 : w^1_i(\lambda) \geq 1\}$.
From (\ref{interw-1w1w}) we obtain $z^1_1(\lambda)\leq z_1(\lambda)$ and $w_i(\lambda)-w^1_i(\lambda)\leq 1$ for $i\geq 1$, hence
$\abs{\bw(\lambda)}-\abs{\bw^1(\lambda)}=\sum_{i=1}^{z_1(\lambda)}(w_i(\lambda)-w_i^1(\lambda))\leq z_1(\lambda)$.
Then,
$\abs{\br^1}-\abs{ \br}+1=\abs{\bW}-\abs{\bW^1}=\sum_{\lambda\in \bar \FF}(\abs{\bw(\lambda)}-\abs{\bw^1(\lambda)})\leq \sum_{\lambda\in \bar \FF}z_1(\lambda)$, i.e.,  (\ref{eqrmins1leqz1}) holds.

Let us prove the converse.
 \begin{enumerate}
 \item 
   Assume that  conditions (\ref{interw-1w1w}) and
   (\ref{eqwr+}) hold.
    From (\ref{interw-1w1w}) we have $\abs{\bW^1}\leq \abs{\bW}$. Then $x=\abs{\bW}-\abs{\bW^1}-1\geq -1$.
    
   \begin{itemize}
   \item If $x=-1$, from  (\ref{eqwr+}) we obtain $\abs{\br}\geq 1$. Then  $c_1=\#\{i\geq 1 : r_1 \geq 1\}\geq 1$. Define
     $$
\begin{array}{ll}
  \widetilde r^1_0 =r_0+1,\\
 \widetilde  r_{c_1}^1=r_{c_1}-1,\\
 \widetilde  r_i^1=r_i, &i \neq 0, c_1.
   \end{array}
$$
 Then $\widetilde r_0^1\geq\widetilde  r_1^1\geq  \dots$. Take  $\widetilde \br^1_*=(\widetilde r_0^1, \widetilde r_1^1, \dots)$.
\item If  $x\geq 0$,
  take
  $\widetilde \br^1_*=\br_*+\overline{(x+1)}$.
     \end{itemize}
   In both cases  $\br_*\cpr \widetilde \br^1_*$
   Define $\widetilde \bs_*^1=\bs_*$. Then
$\abs{\bW^1}+\abs{\widetilde \br^1}+\abs{\widetilde \bs^1}=\abs{\bW^1}+\abs{\br}+x+\abs{\bs}=\abs{\bW}+\abs{\br}+\abs{\bs}-1=\rank(H(s))-1$. We have
$\rank(H(s))-1+\widetilde s_0^1=\rank(H(s))-1+s_0=m-1$ and $\rank(H(s))-1+\widetilde r_0^1=\rank(H(s))+r_0=n$.
   By Lemma \ref{lemmaexistence}
   there exits a pencil $H_1(s)\in \FF[s]^{(m-1)\times n}$ of $\rank(H_1(s))=\rank(H(s))-1$ with Weyr characteristic $(\bW^1, \widetilde \br^1_*, \widetilde  \bs^1_*)$.
   As  condition (\ref{interw-1w1w}), $\br_*\cpr \widetilde \br^1_*$ and  $\bs_*=\widetilde \bs^1_*$ hold, by Proposition \ref{propDox1conj},
there exists a pencil
 $h(s)\in \FF[s]^{1\times n}$ such that $H(s)\se \begin{bmatrix}h(s)\\H_1(s)\end{bmatrix}$.
 \item
  Assume that
  (\ref{colprecconj}) and (\ref{eqrmins1leqz1}) hold.
 From (\ref{eqrmins1leqz1}), $0\leq x= \abs{\br^1}-\abs{\br}+1\leq\sum _{\lambda \in \bar \FF}z_1(\lambda)$.
Let  $x(\lambda)$
be integers such that  $0\leq x(\lambda)\leq z_1(\lambda)$ for  $\lambda\in \bar \FF$, and $\sum _{\lambda \in \bar \FF}x(\lambda)= x$.

Let $\lambda\in \bar \FF$. Define
$$
\begin{array}{ll}
  \widetilde w^1_i(\lambda)=w_i(\lambda), & 1\leq i \leq z_1(\lambda)-x(\lambda),\\
\widetilde   w^1_i(\lambda)=w_i(\lambda)-1, & z_1(\lambda)-x(\lambda)+1\leq i \leq z_1(\lambda),\\
\widetilde   w^1_i(\lambda)=w_i(\lambda)=0, & i>z_1(\lambda).
\end{array}
$$
We have
$\widetilde w^1_1(\lambda)\geq \widetilde w^1_2(\lambda)\geq \dots$. Let $\widetilde \bw^1(\lambda)=(\widetilde w^1_1(\lambda), \dots)$. Then $\Lambda(\widetilde \bw^1)$ is a finite set, $w_i(\lambda)-1\leq  \widetilde w^1_i(\lambda)\leq  w_i(\lambda)$   for $i\geq 1$ and $\lambda \in \bar \FF$, and
$
\abs{\bw(\lambda)} -\abs{\widetilde \bw^1(\lambda)} =x(\lambda);
$
hence 
$$\abs{\bW} -\abs{\widetilde \bW^1} =\sum_{\lambda \in \bar \FF}\abs{\bw(\lambda)}-\sum_{\lambda \in \bar \FF}\abs{\widetilde \bw^1(\lambda)} =\sum_{\lambda \in \bar \FF}x(\lambda)=x=\abs{\br^1}-\abs{\br}+1.$$

Define $\widetilde \bs^1_*=\bs_*$.
Then
$\abs{\widetilde \bW^1}+\abs{\br^1}+\abs{\widetilde \bs^1}=\abs{\bW}+\abs{\br}+\abs{\bs}-1=\rank(H(s))-1.$
The result follows as in item \ref{ittheo11hw}. 
\item
  Assume that  (\ref{roweqconj})
   and
  (\ref{eqabswabsrgeq1}) hold.
  Then
   $\abs{\bW}+\abs{\br}>0$; i.e., $\abs{\bW}>0$ or $\abs{\br}>0$.
  
\begin{itemize}
\item
  If  $\abs{\br}=0$ take  $\widetilde \br^1_*=(r_0+1, 0, \dots)$.
  We have $\abs{\bW}>0$. Let $\lambda_0\in \bar \FF$ be such that $\bw(\lambda_0)\neq 0$ and let $\bz(\lambda_0)=\overline{\bw(\lambda_0)}$.
  Then $z_1(\lambda_0)\geq 1$.
  Define $\widetilde \bw^1(\lambda)=\bw(\lambda)$ for $\lambda \in \bar \FF\setminus \{\lambda_0\}$, and
  $\widetilde \bw^1(\lambda_0)=(\widetilde w^1_1(\lambda_0),\widetilde  w^1_2(\lambda_0), \dots)$ as follows:
  $$
  \begin{array}{ll}
   \widetilde  w^1_i(\lambda_0)=w_i(\lambda_0), & i\neq z_1(\lambda_0),\\
   \widetilde  w^1_i(\lambda_0)=w_i(\lambda_0)-1, & i= z_1(\lambda_0). \\
    \end{array}
  $$
  Then 
    $\Lambda(\widetilde \bw^1)$ is a finite set and
$\abs{\widetilde \bW^1}+\abs{\widetilde \br^1}+\abs{\bs^1}=
\abs{\bW}-1+\abs{\br}+\abs{\bs}=\rank(H(s))-1$.

\item
  If  $\abs{\br}>0$, take 
$\widetilde \bW^1=\bW$.   Let $\bc=\overline{\br}$. Define $\widetilde \br_*^1=(r_0^1, r_1^1, \dots)$ as follows:
  $$
   \begin{array}{ll}
   \widetilde   r^1_0=r_0+1,\\
  \widetilde   r^1_i=r_i, &  i\neq 0,c_1,\\
   \widetilde  r^1_{c_1}=r_{c_1}-1.\\
    \end{array}
  $$
 Then  $\abs{\widetilde \bW^1}+\abs{\widetilde \br^1}+\abs{\bs^1}=
\abs{\bW}+\abs{\br}-1+\abs{\bs}=\rank(H(s))-1$. 
\end{itemize}
\end{enumerate}
In both cases, 
  $w_i(\lambda)-1\leq  \widetilde w^1_i(\lambda)\leq  w_i(\lambda)$   for $i\geq 1$ and $\lambda \in \bar \FF$, and 
 $\br_*\cpr \widetilde \br^1_*$.
The result follows as in item \ref{ittheo11hw}.

\hfill $\Box$

The following theorems (Theorems \ref{theo01h1} and \ref{theo11h1}) can be derived from \cite[Theorem 4.10]{AmBaMaRo24} and \cite[Corollaries 3.14 and 3.16]{AmBaMaRo24_2}.
\begin{theorem}
  \label{theo01h1} 
  Let $H_1(s)\in \FF[s]^{ (m-1)\times n}$ be a matrix pencil with  Weyr characteristic $(\bW^1, \br^1_*,  \bs^1_*)$,
  and for $\lambda \in \bar \FF$ let
$z_1^1(\lambda)=\#\{i\geq 1 : w^1_i(\lambda) \geq 1\}$.
  For $\lambda \in \bar \FF$,
  let  $\bw(\lambda)=(w_1(\lambda), w_2(\lambda), \dots )$ be partitions such that  $\Lambda(\bw)$ is a finite set, $\Lambda(\bw)=\{\lambda_1, \dots, \lambda_\ell\}$,
  and let $\br_*=(r_0, r_1, \dots)$ and $\bs_*=(s_0, s_1, \dots)$  be partitions. 
 There exist pencils 
 $H(s)\in \FF[s]^{m\times n}$, $h(s)\in \FF[s]^{1\times n}$ such that $H(s)\se \begin{bmatrix}h(s)\\H_1(s)\end{bmatrix}$,
$\rank(H(s))=\rank (H_1(s))$, and   
 \begin{enumerate}
 \item \label{ittheo01h1w} 
   $\bW=(\bw(\lambda_1), \dots, \bw(\lambda_\ell))$ is the Weyr characteristic  of the regular part of $H(s)$ if and only if (\ref{interww1w+1}).
   
 \item \label{ittheo01h1r} 
   $\br_*$
   is the column Brunovsky partition of $H(s)$ if and only if (\ref{coleqconj}).

 \item \label{ittheo01h1s} 
   $\bs_*$
   is the row Brunovsky partition of $H(s)$ if and only if (\ref{rowprecconj})
   and
  \begin{equation}\label{eqqsmins1leqz}
     0\leq \abs{\bs}-\abs{ \bs^1}\leq \sum_{\lambda \in \bar \FF}z^1_1(\lambda).
\end{equation}
 
  \end{enumerate}
\end{theorem}
{\bf Proof.} 
Assume that there exist pencils 
$H(s)\in \FF[s]^{m\times n}$, $h(s)\in \FF[s]^{1\times n}$ such that $H(s)\se \begin{bmatrix}h(s)\\H_1(s)\end{bmatrix}$,
$\rank(H(s))=\rank (H_1(s))$,
and   $(\bW, \br_*,  \bs_*)$ is the
Weyr characteristic of $H(s)$. By Proposition \ref{propDox0conj},
(\ref{interww1w+1})--(\ref{rowprecconj}) hold. Since $\rank(H(s))=\rank (H_1(s))$,
we have
$\abs{\bW}+\abs{\br}+\abs{\bs}=\abs{\bW^1}+\abs{\br^1}+\abs{\bs^1}$. From 
(\ref{interww1w+1}) and 
(\ref{coleqconj})
we obtain
$\abs{\bs}-\abs{\bs^1}=\abs{\bW^1}-\abs{\bW}\geq 0$.

For $\lambda \in \bar \FF$,  let $z_1(\lambda)=\#\{i\geq 1 : w_i(\lambda) \geq 1\}$.
From condition (\ref{interww1w+1}) we obtain $z_1(\lambda)\leq z_1^1(\lambda)$ and $w_i^1(\lambda))-w_i(\lambda)\leq 1$, hence
$\abs{\bw^1(\lambda)}-\abs{\bw(\lambda)}=\sum_{i=1}^{z_1^1(\lambda)}(w_i^1(\lambda)-w_i(\lambda))\leq z_1^1(\lambda)$.
Therefore,
 $\abs{\bs}-\abs{ \bs^1}=\abs{\bW^1}-\abs{\bW}=\sum_{\lambda\in \bar \FF}(\abs{\bw^1(\lambda)}-\abs{\bw(\lambda)})\leq \sum_{\lambda\in \bar \FF}z_1^1(\lambda)$, i.e.,  (\ref{eqqsmins1leqz}) holds.

Let us prove the converse. 
   \begin{enumerate}
 \item
   Assume that  (\ref{interww1w+1}) holds. Then $x=\abs{\bW^1}-\abs{\bW}\geq 0$.
    Define $\widetilde \br_*=\br^1_*$ and $\widetilde \bs_*=\bs^1_*+\overline{(x +1)}$. Then 
$\bs^1_*\cpr \widetilde \bs_*$ and
$\abs{\bW}+\abs{\widetilde \br}+\abs{\widetilde \bs}=\abs{\bW^1}-x+\abs{\br^1}+\abs{\bs^1}+x=\rank(H_1(s))$. We have
$\rank(H_1(s))+\widetilde s_0=\rank(H_1(s))+s^1_0+1=m$ and $\rank(H_1(s))+\widetilde r_0=\rank(H_1(s))+r^1_0=n$.
   By Lemma \ref{lemmaexistence}
   there exits a pencil $H(s)\in \FF[s]^{m\times n}$, of $\rank(H(s))=\rank(H_1(s))$ with Weyr characteristic $(\bW, \widetilde \br_*,  \widetilde \bs_*)$.
   As  condition 
 (\ref{interww1w+1}),   $\widetilde \br_*=\br^1_*$ and $\bs^1_*\cpr \widetilde \bs_*$ hold,
   by Proposition \ref{propDox0conj},
there exists a pencil
 $h(s)\in \FF[s]^{1\times n}$ such that $H(s)\se \begin{bmatrix}h(s)\\H_1(s)\end{bmatrix}$.
\item
  Assume that (\ref{coleqconj}) holds.
Define $\widetilde \bs_*=(\widetilde s_0, \widetilde s_1, \dots)$ as
$$
\begin{array}{ll}
\widetilde s_0=s^1_0+1,\\
\widetilde s_i=s^1_i, &i\geq 1.\\
\end{array}
$$
Then $\bs^1_*\cpr \widetilde \bs_*$ and  $\abs{\widetilde \bs}=\abs{\bs^1}$. Define $\widetilde \bW=\bW^1$. 
The result follows as in item \ref{ittheo01h1w}.
\item 
  Assume that (\ref{rowprecconj})
   and
   (\ref{eqqsmins1leqz}) hold.
 From (\ref{eqqsmins1leqz}), $0\leq y =\abs{\bs}-\abs{\bs^1}\leq\sum _{\lambda \in \bar \FF}z_1^1(\lambda)$.
Let  $y(\lambda)$
be integers such that  $0\leq y(\lambda)\leq z_1^1(\lambda)$, for  $\lambda\in \bar \FF$, and $\sum _{\lambda \in \bar \FF}y(\lambda)= y$.

For $\lambda\in \bar \FF$ define
$$
\begin{array}{ll}
 \widetilde  w_i(\lambda)=w_i^1(\lambda), & 1\leq i \leq z_1^1(\lambda)-y(\lambda),\\
\widetilde   w_i(\lambda)=w_i^1(\lambda)-1, & z_1^1(\lambda)-y(\lambda)+1\leq i \leq z_1^1(\lambda),\\
 \widetilde  w_i(\lambda)=w_i^1(\lambda)=0, & i>z_1^1(\lambda).
\end{array}
$$
Then
$\widetilde w_1(\lambda)\geq \widetilde w_2(\lambda)\geq \dots$. Let $\widetilde \bw(\lambda)=(\widetilde w_1(\lambda), \dots)$. Then $\Lambda(\widetilde \bw)$ is a finite set,  $\widetilde w_i(\lambda)\leq  w^1_i(\lambda)\leq \widetilde w_i(\lambda)+1$ for  $i\geq 1$ and $ \lambda \in \bar \FF$, and
$
\abs{\bw^1(\lambda)} -\abs{\widetilde \bw(\lambda)} =y(\lambda);
$
hence 
$\abs{\bW^1} -\abs{\widetilde\bW} =\sum_{\lambda \in \bar \FF}\abs{\bw^1(\lambda)}-\sum_{\lambda \in \bar \FF}\abs{\bw(\lambda)} =\sum_{\lambda \in \bar \FF}y(\lambda)=y=\abs{\bs}-\abs{\bs^1}.$
Define $\widetilde \br_*=\br_*^1$. Then 
$\abs{\widetilde\bW}+\abs{\widetilde\br}+\abs{\bs}=\abs{\bW^1}+\abs{\br^1}+\abs{\bs^1}.$
The result follows as in item \ref{ittheo01h1w}.
  \end{enumerate}
\hfill $\Box$

\begin{theorem}
  \label{theo11h1} 
  Let $H_1(s)\in \FF[s]^{ (m-1)\times n}$ be a matrix pencil with  Weyr characteristic $(\bW^1, \br^1_*,  \bs^1_*)$, and let  $\bc^1 =(c^1_1, \dots )$ be  the conjugate partition of $\br^1=(r^1_1, \dots )$.
  For $\lambda \in \bar \FF$,
  let  $\bw(\lambda)=(w_1(\lambda), w_2(\lambda), \dots )$ be partitions such that  $\Lambda(\bw)$ is a finite set, $\Lambda(\bw)=\{\lambda_1, \dots, \lambda_\ell\}$,
  and let $\br_*=(r_0, r_1, \dots)$ and $\bs_*=(s_0, s_1, \dots)$  be partitions. 
 There exist pencils 
 $H(s)\in \FF[s]^{m\times n}$, $h(s)\in \FF[s]^{1\times n}$ such that $H(s)\se \begin{bmatrix}h(s)\\H_1(s)\end{bmatrix}$,
$\rank(H(s))=\rank (H_1(s))+1$, and   
 \begin{enumerate}
 \item \label{ittheo11h1w}
$\bW=(\bw(\lambda_1), \dots, \bw(\lambda_\ell))$ is the Weyr characteristic  of the regular part of $H(s)$ if and only if (\ref{interw-1w1w}),
\begin{equation}\label{eqr10}
     r_0^1\geq 1,
   \end{equation}
\begin{equation}\label{eqwr1}
\mbox{if $r_0^1=1$, then }
     \abs{\bW}-\abs{\bW^1}-1=c^1_1, 
   \end{equation}
   and
 \begin{equation}\label{eqwr}  
 \mbox{if $r_0^1>1$, then }  
     \abs{\bW}-\abs{\bW^1}-1=c^1_1 \mbox{ or } \abs{\bW}-\abs{\bW^1}-1\leq c^1_2.
   \end{equation}
 \item \label{ittheo11h1r} 
   $\br_*$
   is the column Brunovsky partition of $H(s)$ if and only if (\ref{colprecconj})
and  \begin{equation}\label{eqrr1}
     0\leq\abs{\br^1} -\abs{\br}+1.
   \end{equation}
 \item \label{ittheo11h1s} 
   $\bs_*$
   is the row Brunovsky partition of $H(s)$ if and only if (\ref{roweqconj})
   and 
(\ref{eqr10}).

  \end{enumerate}
\end{theorem}
{\bf Proof.}
Assume that there exist pencils 
$H(s)\in \FF[s]^{m\times n}$, $h(s)\in \FF[s]^{1\times n}$ such that $H(s)\se \begin{bmatrix}h(s)\\H_1(s)\end{bmatrix}$,
$\rank(H(s))=\rank (H_1(s))+1$,
and   $(\bW, \br_*,  \bs_*)$ is the
Weyr characteristic of $H(s)$. By Proposition \ref{propDox1conj} and Remark \ref{remranks}, 
(\ref{interw-1w1w})--(\ref{roweqconj}) and (\ref{eqrr1}) hold. By (\ref{roweqconj}),  since $\rank(H(s))=\rank (H_1(s))+1$,
we have $\abs{\br}=\abs{\br^1}+\abs{\bW^1}-\abs{\bW}+1$. Then 
(\ref{eqr10})--(\ref{eqwr}) follow from (\ref{colprecconj}) by Lemma \ref{lemmatec1}.

Let us prove the converse.
 \begin{enumerate}
 \item
   Assume that  (\ref{interw-1w1w}),  and
(\ref{eqr10})--(\ref{eqwr})  hold.
   By Lemma \ref{lemmatec1},
   there exists a partition $\widetilde \br_*=(\widetilde r_0, \widetilde r_1, \dots )$ such that $\widetilde \br_*\cpr\br^1_*$ and $\abs{\widetilde \br}=\abs{\br^1}-\abs{\bW}+ \abs{\bW^1}+1$.
 Define $\widetilde \bs_*=\bs^1_*$. 
   Then $\abs{\bW}+\abs{\widetilde \br}+\abs{\widetilde \bs}=\abs{\bW^1}+\abs{\br^1}+\abs{\bs^1}+1=\rank(H_1(s))+1$. We have
$\rank(H_1(s))+1+\widetilde s_0=\rank(H_1(s))+s^1_0+1=m$ and $\rank(H_1(s))+1+\widetilde r_0=\rank(H_1(s))+r^1_0=n$.
   By Lemma \ref{lemmaexistence}
   there exits a pencil $H(s)\in \FF[s]^{m\times n}$, of $\rank(H(s))=\rank(H_1(s))+1$ with Weyr characteristic $(\bW, \widetilde\br_*, \widetilde \bs_*)$.
   As condition
   (\ref{interw-1w1w}), $\widetilde \br_*\cpr \br^1_*$ and $\widetilde \bs_*=\bs^1_*$ hold,
   by Proposition \ref{propDox1conj},
   there exist a pencil $h(s)\in \FF[s]^{1\times n}$ such that $H(s)\se \begin{bmatrix}h(s)\\H_1(s)\end{bmatrix}$.
 \item Assume that (\ref{colprecconj})
and  (\ref{eqrr1}) hold.
Then $x= \abs{\br^1}- \abs{\br}+1\geq 0$.
 Fix $\lambda_0\in \bar \FF$ and for $\lambda \in \bar \FF$ define a partition $\widetilde \bw(\lambda)=(\widetilde w_1(\lambda),\widetilde  w_2(\lambda), \dots)$ as follows:
$$
\begin{array}{ll}
  \widetilde \bw(\lambda) =\bw^1(\lambda), \quad \lambda\in \bar \FF\setminus\{\lambda_0\},\\
 \widetilde  \bw(\lambda_0)=\bw^1(\lambda_0)+\overline{(x )},
  \end{array}
$$
 and define $\widetilde \bs_*=\bs^1_*$. 
The result follows as in item \ref{ittheo11hw}. 
 \item Assume that  (\ref{roweqconj})
   and 
   (\ref{eqr10}) hold.
  Define $\widetilde \br_*=(\widetilde r_0, \widetilde r_1, \dots)$ as follows:
$$
\begin{array}{ll}
\widetilde r_i=r^1_i-1, &0\leq i \leq c^1_1,\\
\widetilde r_i=r^1_i=0, &i> c^1_1.\\
\end{array}
$$
Then $\widetilde \br_*\cpr \br^1$ and  $\abs{\widetilde \br}=\abs{\br^1}-c^1_1$. 

Fix $\lambda_0\in \bar \FF$ and for $\lambda \in \bar \FF$ define a partition $\widetilde \bw(\lambda)=(\widetilde w_1(\lambda),\widetilde  w_2(\lambda), \dots)$ as follows:
$$
\begin{array}{ll}
\widetilde   \bw(\lambda) =\bw^1(\lambda), \quad \lambda\in \bar \FF\setminus\{\lambda_0\},\\
\widetilde   \bw(\lambda_0)=\bw^1(\lambda_0)+\overline{(c^1_1 +1)},
  \end{array}
$$
The result follows as in item \ref{ittheo11hw}. 

 \end{enumerate}
 \hfill $\Box$

\section{
Reachability of the bounds
}
\label{reachability} 

In this section we show with different examples that all the bounds obtained  in Theorem \ref{theogpertboundsrow}  are reachable. We are assuming the hypothesis of Theorem \ref{theogpertboundsrow}.

\begin{enumerate}
\item  If $\rank( H(s)+P(s))=\rank( H(s))$, then (\ref{eqgpboundwrrr1}), (\ref{eqgpboundrrrr1+})  and (\ref{eqgpboundsrrr1}) hold.    
  \begin{itemize}
  \item \underline{Bounds in (\ref{eqgpboundwrrr1}).}
  Let $H(s)\in \FF[s]^{m\times n}$ be  a pencil  with Weyr characteristic $(\bW, \br_*, \bs_*)$, where $\bw(0)=(1,1,0)$, $\bw(\lambda)=(0)$ for $\lambda \neq 0$, and $r_0\geq 1$.  
     
    Let $w^1(0)=(1,0)$ and $w^1(\lambda)=(0)$ for $\lambda \neq 0$. 
  Then   (\ref{interw-1w1w}) holds and, as $\abs{\bW^1}-\abs{\bW}+1=0$, 
   (\ref{eqwr+}) also holds.
   By Theorem \ref{theo011h} 
  there exist pencils 
 $H_1(s)\in \FF[s]^{(m-1)\times n}$, $h(s)\in \FF[s]^{1\times n}$ such that $H(s)\se \begin{bmatrix}h(s)\\H_1(s)\end{bmatrix}$,
$\rank(H_1(s))=\rank(H(s))-1$, and   
 $\bW^1$ is the Weyr characteristic of the regular part of $H_1(s)$. Observe that 
 if  $\br^1_*$ is the column Brunovsky partition of $H_1(s)$, then
 $r_0^1=n-\rank(H_1(s))=n-\rank(H(s))+1=r_0+1\geq 2$.

Additionally, take $\widehat \bW(0)=(2,0)$ and $\widehat \bW(\lambda)=(0)$ for $\lambda \neq 0$.  We have
$\widehat w_i(\lambda)-1\leq  w^1_i(\lambda)\leq  \widehat w_i(\lambda)$ for $i\geq 1$ and $\lambda \in \bar \FF$. Let $\bc^1=(c_1^1, \dots)$ be the conjugate partition of $\br^1$. Then $\abs{ \widehat  \bW}-\abs{\bW^1}-1=0\leq c_2^1$.
    By Theorem \ref{theo11h1} 
 there exist pencils 
 $G(s)\in \FF[s]^{m\times n}$, $g(s)\in \FF[s]^{1\times n}$ such that $G(s)\se \begin{bmatrix}g(s)\\H_1(s)\end{bmatrix}$,
$\rank(G(s))=\rank (H_1(s))+1$, and
$\widehat \bW$ is the Weyr characteristic of the regular part of $G(s)$.

By Lemma \ref{lemmaeq}
there exists a matrix pencil  $P(s)=uv(s)^T \in \FF[s]^{m\times n}$ such that
$H(s)+P(s)\se G(s)$; therefore, $\widehat \bW$ is the Weyr characteristic of the regular part of $H(s)+P(s)$.
We have $\rank(H(s)+P(s))=\rank(G(s))=\rank (H_1(s))+1=\rank( H(s))$,
$\widehat w_1(0)-w_1(0)=1$ and $\widehat w_2(0)-w_2(0)=-1$; i.e. the bounds in (\ref{eqgpboundwrrr1}) are reachable.

Furthermore, assume that $\bs_*=(1,1,0)$ and take  $\bw^1(0)=(2,1,0)$ and $\bw^1(\lambda)=(0)$ for $\lambda \neq 0$. Applying Theorems \ref{theo01h} and  \ref{theo01h1} we obtain a similar result with 
$\rank(H(s)+P(s))=\rank (H(s))= \rank (H_1(s))$.

\item \underline{Bounds in (\ref{eqgpboundrrrr1+}) when $\br= (0)$.}
 Let $H(s)\in \FF[s]^{m\times n}$ be  a pencil  with Weyr characteristic $(\bW, \br_*, \bs_*)$, where
$\br_*=(2,0)$,  $\bw(0)=(1,1,0)$ and  $\bw(\lambda)=(0)$ for $\lambda \neq 0$. Take $\br^1_*=(3,1, 0)$.
 Then   $\abs{\br^1}=1$, $\abs{\br}=0$,  and (\ref{colprecconj}) and
 (\ref{eqrmins1leqz1}) hold.
   By Theorem \ref{theo011h}
  there exist pencils 
 $H_1(s)\in \FF[s]^{(m-1)\times n}$, $h(s)\in \FF[s]^{1\times n}$ such that $H(s)\se \begin{bmatrix}h(s)\\H_1(s)\end{bmatrix}$,
$\rank(H_1(s))=\rank (H(s))-1$, and   
   $\br_*^1$
   is the column Brunovsky partition of $H_1(s)$.

Moreover, let  $\widehat \br_*=(2,2,0)$. Then $\abs{\widehat \br}=2$, $\widehat \br_*\cpr \br_*^1$ and $\abs{ \br^1}-\abs{ \widehat \br}+1=0$.
As in the previous item, by Theorem \ref{theo11h1} and Lemma \ref{lemmaeq},
there exists  a  pencil  $P(s)=uv(s)^T \in \FF[s]^{m\times n}$ such that
 $\rank(H(s)+P(s))=\rank (H_1(s))+1=\rank( H(s))$ and $\widehat \br_*$
   is the column Brunovsky partition of  $H(s)+P(s)$. As
$\widehat r_{1}-r_{1}=2$ and $\widehat r_{2}-r_{2}=0$,  the  bounds in (\ref{eqgpboundrrrr1+}) are reachable when $\br=(0)$.

\item \underline{Lower bound in (\ref{eqgpboundrrrr1+}) when $\br\neq  (0)$.}
Let $H(s)\in \FF[s]^{m\times n}$ be  a pencil  with Weyr characteristic $(\bW, \br_*, \bs_*)$, where
$\br_*=(\overbrace{11, \dots, 11,}^{11}1,1,1, 0)$, and let 
$$
\br^1_*=(\overbrace{12, \dots, 12,}^{10}1,1,1,1, 0), \quad \widehat \br_*=(\overbrace{11, \dots, 11,}^{10} 0).$$
With a similar reasoning as in the previous item, applying Theorems \ref{theo011h} and
  \ref{theo11h1}  we obtain that there is a pencil  $P(s)=uv(s)^T \in \FF[s]^{m\times n}$ such that
$\rank(H(s)+P(s))=\rank( H(s))$, and $\br_*$ and $\widehat \br_*$ are the column Brunovsky partitions of $H(s)$ and $H(s)+P(s)$, respectively.

As
$\widehat  r_{10}-r_{10}=0-11= -11=- \left [ \sqrt{113}\right ]-1=-\left [ \sqrt{\abs{\br}}\right ]-1$, when $\br\neq (0)$  the lower bound in (\ref{eqgpboundrrrr1+}) is reachable.

\item \underline{Upper bound in (\ref{eqgpboundrrrr1+}) when $\br\neq  (0)$.} 
Let $H(s)\in \FF[s]^{m\times n}$ be  a pencil  with Weyr characteristic $(\bW, \br_*, \bs_*)$, where
$$
\br_*=( \overbrace{10, \dots, 10,}^{9} 0);  \quad \bw(0)=\overline{(10,0)},\quad  \bw(\lambda)=(0)\mbox{ for $\lambda \neq 0$},$$
and let 
$$
\br^1_*=(\overbrace{11, \dots, 11,}^{9}1, 0), \quad \widehat \br_*=(\overbrace{10, \dots, 10,}^{10}0).$$
Proceeding as in the previous items, we obtain that there is a  pencil $P(s)=uv(s)^T \in \FF[s]^{m\times n}$ such that
$\rank(H(s)+P(s))=\rank( H(s))$, and $\br_*$ and $\widehat \br_*$ are the column Brunovsky partitions of $H(s)$ and $H(s)+P(s)$, respectively.

Here,  
$\widehat r_{9}-r_{9}=10-0=10=  \left [ \sqrt{79}\right ]+2=\left [ \sqrt{\abs{\br}-1}\right ]+2$;  i.e. the upper bound in (\ref{eqgpboundrrrr1+}) is reachable when $\br\neq (0)$.

\item \underline{Bounds in (\ref{eqgpboundsrrr1}) when $\bs= (0)$.}
 Let $H(s)\in \FF[s]^{m\times n}$ be  a pencil  with Weyr characteristic $(\bW, \br_*, \bs_*)$, where
$\bs_*=(2,0)$,  $\bw(0)=(1,0)$ and  $\bw(\lambda)=(0)$ for $\lambda \neq 0$. Take $\bs^1_*=(1, 0)$.
 Then    (\ref{rowprecconj}) and
 (\ref{abss1leqabss}) hold.
   By Theorem \ref{theo01h}
  there exist pencils 
 $H_1(s)\in \FF[s]^{(m-1)\times n}$, $h(s)\in \FF[s]^{1\times n}$ such that $H(s)\se \begin{bmatrix}h(s)\\H_1(s)\end{bmatrix}$,
$\rank(H_1(s))=\rank (H(s))$, and   
   $\bs_*^1$
   is the row Brunovsky partition of $H_1(s)$.
   Observe that if $\bW^1=(\bw^1(\lambda_1), \dots, \bw^1(\lambda_\ell))$ is the Weyr characteristic of the regular part of $H_1(s)$, then, since $\rank(H_1(s))=\rank (H(s))$, by Proposition \ref{propDox0conj}    we obtain 
   $\abs{\bW^1}=\abs{\bW}+\abs{\bs}-\abs{\bs^1}=1$ and $w_1(0)=1\leq w^1(0)$. Hence, $\bw^1(0)=(1,0)$ and  $\bw^1(\lambda)=(0)$ for $\lambda \neq 0$. 

Let  $\widehat \bs_*=(2,1,0)$. Then $\bs^1_*\cpr \widehat \bs_*$ and $\abs{ \widehat \bs}-\abs{ \bs^1}=1=\sum_{\lambda \in \bar \FF}z^1_1(\lambda) $, where
$z_1^1(\lambda)=\#\{i\geq 1 : w^1_i(\lambda) \geq 1\}$ for 
   $\lambda \in \bar \FF$.  
 By Theorem \ref{theo01h1}
 there exist pencils 
 $G(s)\in \FF[s]^{m\times n}$, $g(s)\in \FF[s]^{1\times n}$ such that $G(s)\se \begin{bmatrix}g(s)\\H_1(s)\end{bmatrix}$,
$\rank(G(s))=\rank (H_1(s))$, and
$\widehat \bs_*$
   is the row Brunovsky partition of $G(s)$.
    By Lemma \ref{lemmaeq}
there exists a matrix pencil  $P(s)=uv(s)^T \in \FF[s]^{m\times n}$ such that
$H(s)+P(s)\se G(s)$; it means that $\widehat \bs_*$
   is the row Brunovsky partition of  $H(s)+P(s)$.
As $\rank(H(s)+P(s))=\rank(G(s))=\rank (H_1(s))=\rank( H(s))$,
$\widehat s_{0}-s_{0}=0$ and $\widehat s_{1}-s_{1}=1$, the  bounds in (\ref{eqgpboundsrrr1}) are reachable when $\bs=(0)$.

\item \underline{Lower bound in  (\ref{eqgpboundsrrr1}) when $\bs\neq (0)$.}
 Let $H(s)\in \FF[s]^{m\times n}$ be  a pencil  with Weyr characteristic $(\bW, \br_*, \bs_*)$, where
$\bs_*=(\overbrace{101, \dots, 101,}^{102} 0)$. Let $\bs^1_*=(\overbrace{100, \dots, 100,}^{102} 0)$.
 Then,   $\abs{\bs_*}=10302$, $\abs{\bs^1_*}=10200$,  and 
   (\ref{rowprecconj}) and 
   (\ref{abss1leqabss}) hold. 
   By Theorem \ref{theo01h}
  there exist pencils 
 $H_1(s)\in \FF[s]^{(m-1)\times n}$, $h(s)\in \FF[s]^{1\times n}$ such that $H(s)\se \begin{bmatrix}h(s)\\H_1(s)\end{bmatrix}$,
$\rank(H_1(s))=\rank (H(s))$, and   
   $\bs_*^1$
   is the row Brunovsky partition of $H_1(s)$.

Take  $\widehat \bs_*=(\overbrace{101, \dots, 101,}^{101} 0)$. Then, $\bs^1_*\cpr \widehat \bs_*$  and $\abs{\widehat \bs}-\abs{ \bs^1}=10100-10100=0$.
As in the previous item there exists a matrix pencil  $P(s)=uv(s)^T \in \FF[s]^{m\times n}$ such that
$\rank(H(s)+P(s))=\rank( H(s))$,
and $\widehat \bs_*$
   is the row Brunovsky partition of  $H(s)+P(s)$.
As
$\widehat s_{101}-s_{101}=0-101=-101=-\left [ \sqrt{10302}\right ]=-\left [ \sqrt{\abs{\bs_*}}\right ]$,  the  lower bound in (\ref{eqgpboundsrrr1}) is reachable.

\item \underline{ Upper bound in (\ref{eqgpboundsrrr1}) when $\bs\neq  (0)$.}
 Let $H(s)\in \FF[s]^{m\times n}$ be  a pencil  with Weyr characteristic $(\bW, \br_*, \bs_*)$, where
$\bs_*=(\overbrace{10, \dots, 10,}^{10} 0)$,  $\bw(0)=\overline{(10,0)}=(\overbrace{1, \dots, 1,}^{10} 0)$ and  $\bw(\lambda)=(0)$ for $\lambda \neq 0$. Take $\bs^1_*=(\overbrace{9, \dots, 9,}^{11} 0)$.
 Then    (\ref{rowprecconj}) and
 (\ref{abss1leqabss}) hold.
   By Theorem \ref{theo01h}
  there exist pencils 
 $H_1(s)\in \FF[s]^{(m-1)\times n}$, $h(s)\in \FF[s]^{1\times n}$ such that $H(s)\se \begin{bmatrix}h(s)\\H_1(s)\end{bmatrix}$,
$\rank(H_1(s))=\rank (H(s))$, and   
   $\bs_*^1$
   is the row Brunovsky partition of $H_1(s)$.
   If $\bW^1$ is the Weyr characteristic of the regular part of $H_1(s)$, then, since $\rank(H_1(s))=\rank (H(s))$, by Proposition \ref{propDox0conj} we obtain 
   $\abs{\bW^1}=\abs{\bW}+\abs{\bs}-\abs{\bs^1}=\abs{\bW}=10$ and $w_i(0)=1\leq w^1_i(0)$ for $1\leq i\leq 10$, hence 
   $\bw^1(0)=\overline{(10,0)}$ and  $\bw^1(\lambda)=(0)$ for $\lambda \neq 0$; it means that $z_1^1(0)=\#\{i\geq 1 : w^1_i(0) \geq 1\}=10$.

Let  $\widehat \bs_*=(\overbrace{10, \dots, 10,}^{11} 0)$. Then $\bs^1_*\cpr \widehat \bs_*$ and $\abs{ \widehat \bs}-\abs{ \bs^1}=10$.
As in the previous items, 
 by Theorem \ref{theo01h1} and 
 Lemma \ref{lemmaeq}
there exists a matrix pencil  $P(s)=uv(s)^T \in \FF[s]^{m\times n}$ such that
$\rank(H(s)+P(s))=\rank( H(s))$ and 
  $\widehat \bs_*$
   is the row Brunovsky partition of  $H(s)+P(s)$.
As 
$\widehat s_{10}-s_{10}=10-0=10=\left [ \sqrt{100}\right ]=\left [ \sqrt{\abs{\bs_*}}\right ]$,  the  upper bound in (\ref{eqgpboundsrrr1}) is reachable.

\end{itemize}

\item If $\rank (H(s)+P(s))=\rank ( H(s))-1$, then 
    (\ref{eqgpboundwrr+})--(\ref{eqgpboundsrr+}) hold.    
\begin{itemize}
\item \underline{Bounds in (\ref{eqgpboundwrr+}).}
  Let $H(s)\in \FF[s]^{m\times n}$ be  a pencil  with Weyr characteristic $(\bW, \br_*, \bs_*)$, where $\bw(0)=(2,0)$  and $\bw(\lambda)=(0)$ for $\lambda \neq 0$.
     
    Let $\bw^1(0)=(1,0)$  and $\bw^1(\lambda)=(0)$ for $\lambda \neq 0$.
  Then   (\ref{interw-1w1w}) holds and, as $\abs{\bW^1}-\abs{\bW}+1=0$, 
   (\ref{eqwr+}) also holds.
   By Theorem \ref{theo011h} 
  there exist pencils 
 $H_1(s)\in \FF[s]^{(m-1)\times n}$, $h(s)\in \FF[s]^{1\times n}$ such that $H(s)\se \begin{bmatrix}h(s)\\H_1(s)\end{bmatrix}$,
$\rank(H_1(s))=\rank (H(s))-1$, and   
 $\bW^1$ is the Weyr characteristic of the regular part of $H_1(s)$.

Take $\widehat \bW(\lambda)=(0)$ for $\lambda \in \bar \FF$.
Then 
$\widehat w_i(\lambda)\leq  w^1_i(\lambda)\leq  \widehat w_i(\lambda)+1$ for $i\geq 1$ and $\lambda \in \bar \FF$.
    By Theorem \ref{theo01h1}
 there exist pencils 
 $G(s)\in \FF[s]^{m\times n}$, $g(s)\in \FF[s]^{1\times n}$ such that $G(s)\se \begin{bmatrix}g(s)\\H_1(s)\end{bmatrix}$,
$\rank(G(s))=\rank (H_1(s))$, and
$\widehat \bW$ is the Weyr characteristic of the regular part of $G(s)$.

By Lemma \ref{lemmaeq}
there exists a matrix pencil  $P(s)=uv(s)^T \in \FF[s]^{m\times n}$ such that
$H(s)+P(s)\se G(s)$; therefore $\widehat \bW$ is the Weyr characteristic of the regular part of $H(s)+P(s)$.
Observe that $\rank(H(s)+P(s))=\rank(G(s))=\rank (H_1(s))=\rank( H(s))-1$,
$\widehat w_1(0)-w_1(0)=-2$ and $\widehat w_2(0)-w_2(0)=0$; i.e. the bounds in (\ref{eqgpboundwrr+}) are reachable.

\item \underline{Bounds in (\ref{eqgpboundrrr+}).}
Let $H(s)\in \FF[s]^{m\times n}$ be  a pencil  with Weyr characteristic $(\bW, \br_*, \bs_*)$,  where
$
\br_*=( \overbrace{100, \dots, 100,}^{101} 0)$.
Take
$$
\br^1_*=(\overbrace{101, \dots, 101,}^{100} 0), \quad \widehat \br_*=\br^1_*.$$
As in the previous item,
applying Theorems \ref{theo011h} 
  and  \ref{theo01h1}
we obtain that there exists a pencil $P(s)=uv(s)^T \in \FF[s]^{m\times n}$, such that
$\rank(H(s)+P(s))=\rank( H(s))-1$, and $\br_*$ and $\widehat \br_*$ are the column Brunovsky partitions of $H(s)$ and $H(s)+P(s)$, respectively.

Since $\widehat r_{1}-r_{1}=101-100=1$ and 
$\widehat r_{100}-r_{100}=0-100=-100=  -\left [ \sqrt{10000}\right ]=-\left [ \sqrt{\abs{\br}}\right ]$, the bounds in (\ref{eqgpboundrrr+}) are reachable.

\item \underline{Bounds in (\ref{eqgpboundsrr+}).}
Let $H(s)\in \FF[s]^{m\times n}$ be  a pencil  with Weyr characteristic $(\bW, \br_*, \bs_*)$,  where
$
\br_*=(1,1,0)$ and  $\bs_*=(\overbrace{50, \dots, 50,}^{54} 1,0)$. Take  $$\bs^1_*=\bs_*,\quad 
\widehat \bs_*=(\overbrace{51, \dots, 51,}^{52} 50,0).$$
Then $\abs{\bW}+\abs{\br_*}\geq \abs{\br_*}=2>1=r_0$, $\bs_*^1\cpr \widehat \bs_*$ and $\abs{\widehat \bs}-\abs{\bs^1}=2651-2651=0$. 
As in the previous item,
applying Theorems \ref{theo011h} 
  and  \ref{theo01h1}
we obtain that there is a pencil $P(s)=uv(s)^T \in \FF[s]^{m\times n}$ such that
$\rank(H(s)+P(s))=\rank( H(s))-1$, and $\bs_*$ and $\widehat \bs*$ are the row Brunovsky partitions of $H(s)$ and $H(s)+P(s)$, respectively.

Since $\widehat s_{1}-s_{1}=51-50=1$ and 
$\widehat s_{53}-s_{53}=0-50=-50=  1-\left [ \sqrt{2702}\right ]=1-\left [ \sqrt{\abs{\bs_*}+1}\right ]$, the bounds in (\ref{eqgpboundsrr+}) are reachable.

\end{itemize}
\item   If $\rank ( H(s)+P(s))=\rank ( H(s))+1$, then (\ref{eqgpboundwr+r})--
      (\ref{eqgpboundsr+r}) hold.
\begin{itemize}
\item \underline{Bounds in (\ref{eqgpboundwr+r}).}
 Let $H(s)\in \FF[s]^{m\times n}$ be  a pencil  with Weyr characteristic $(\bW, \br_*, \bs_*)$, where  $\bw(\lambda)=(0)$ for $\lambda \in \bar \FF$,   $\br_*=(2,0)$  and  $\bs_*=(1,1,0)$ (therefore, 
 $\bu=\overline{\bs}=(1, 0)$.
     
    Let $\bw^1(0)=(1,0)$  and  $\bw^1(\lambda)=(0)$ for $\lambda \neq 0$.
  Then   (\ref{interww1w+1}) holds, $s_0=1$ and $\abs{\bW^1}-\abs{\bW}=1=u_1$.
   By Theorem \ref{theo01h}
  there exist pencils 
 $H_1(s)\in \FF[s]^{(m-1)\times n}$, $h(s)\in \FF[s]^{1\times n}$ such that $H(s)\se \begin{bmatrix}h(s)\\H_1(s)\end{bmatrix}$,
$\rank(H_1(s))=\rank (H(s))$, and   
 $\bW^1$ is the Weyr characteristic of the regular part of $H_1(s)$. Let $\br^1_*$ be the column Brunovsky partition of $H_1(s)$, then
 $r_0^1=n-\rank(H_1(s))=n-\rank(H(s))=r_0=2$.

Let $\widehat \bw(0)=(2,0)$  and  $\widehat \bw(\lambda)=(0)$ for $\lambda \neq 0$.
Then $\widehat w_i(\lambda)-1\leq  w^1_i(\lambda)\leq \widehat  w_i(\lambda)$ for  $i\geq 1$ and $\lambda \in \bar \FF$,
and  $\abs{\widehat \bW}-\abs{\bW^1}-1=0$.
    By Theorem \ref{theo11h1}
 there exist pencils 
 $G(s)\in \FF[s]^{m\times n}$, $g(s)\in \FF[s]^{1\times n}$ such that $G(s)\se \begin{bmatrix}g(s)\\H_1(s)\end{bmatrix}$,
$\rank(G(s))=\rank (H_1(s))+1$, and
$\widehat \bW$ is the Weyr characteristic of the regular part of $G(s)$.

By Lemma \ref{lemmaeq}
there exists a matrix pencil  $P(s)=uv(s)^T \in \FF[s]^{m\times n}$ such that
$H(s)+P(s)\se G(s)$; therefore $\widehat \bW$ is the Weyr characteristic of the regular part of $H(s)+P(s)$.
As $\rank(H(s)+P(s))=\rank(G(s))=\rank (H_1(s))+1=\rank( H(s))+1$,
$\widehat w_1(0)-w_1(0)=2$ and $\widehat w_2(0)-w_2(0)=0$,  the bounds in (\ref{eqgpboundwr+r}) are reachable.

\item \underline{Bounds in (\ref{eqgpboundrr+r}).}
Let $H(s)\in \FF[s]^{m\times n}$ be  a pencil  with Weyr characteristic $(\bW, \br_*, \bs_*)$,  where
$\br_*=( \overbrace{10, \dots, 10,}^{10} 0)$ and $ \bs_*=(1,0) $. Take 
$$
\br^1_*=\br_*, \quad 
\widehat \br_*=(\overbrace{9, \dots, 9,}^{11} 0).$$
Then $s_0=1$, $\widehat \br_*\cpr \br_*^1$,  and
$\abs{\br^1}-\abs{\widehat \br}+1
=90-90+1=1\geq 0$.

As in the previous item,
applying Theorems \ref{theo01h} 
  and  \ref{theo11h1}
we obtain that there exists a pencil  $P(s)=uv(s)^T \in \FF[s]^{m\times n}$ such that
$\rank(H(s)+P(s))=\rank( H(s))+1$, and $\br_*$ and $\widehat \br_*$ are the column Brunovsky partitions of $H(s)$ and $H(s)+P(s)$, respectively.

We have $\widehat r_{1}-r_{1}=9-10=-1$ and 
$\widehat r_{10}-r_{10}=9-0=9=\left [ \sqrt{91}\right ]=\left [ \sqrt{\abs{\br}+1}\right ]$; i.e. the bounds in (\ref{eqgpboundrr+r}) are reachable.

\item \underline{Bounds in (\ref{eqgpboundsr+r}).}
Let $H(s)\in \FF[s]^{m\times n}$ be  a pencil  with Weyr characteristic $(\bW, \br_*, \bs_*)$,  where
$
\bs_*=( \overbrace{100, \dots, 100,}^{100} 0)$ and $\br_*=(1,0)$. Take $$
\bs^1_*=(\overbrace{99, \dots, 99,}^{101} 0), \quad \widehat \bs_*= \bs^1.$$
Then  (\ref{rowprecconj})
   and
   (\ref{abss1leqabss}) hold.
By Theorem \ref{theo01h} there exist pencils 
 $H_1(s)\in \FF[s]^{(m-1)\times n}$, $h(s)\in \FF[s]^{1\times n}$ such that $H(s)\se \begin{bmatrix}h(s)\\H_1(s)\end{bmatrix}$,
$\rank(H_1(s))=\rank (H(s))$, and   
 $\bs_*^1$ is the row Brunovsky partition of $H_1(s)$.
 Let $\br^1_*$ be the column Brunovsky partition of $H_1(s)$, then
 $r_0^1=n-\rank(H_1(s))=n-\rank(H(s))=r_0=1$.

By Theorem \ref{theo11h1},
 there exist pencils 
 $G(s)\in \FF[s]^{m\times n}$, $g(s)\in \FF[s]^{1\times n}$ such that $G(s)\se \begin{bmatrix}g(s)\\H_1(s)\end{bmatrix}$,
$\rank(G(s))=\rank (H_1(s))+1$, and
  $\widehat \bs_*$ is the row Brunovsky partition of  of $G(s)$.

By Lemma \ref{lemmaeq}
there exists a matrix pencil  $P(s)=uv(s)^T \in \FF[s]^{m\times n}$ such that
$H(s)+P(s)\se G(s)$; therefore $\widehat \bs_*$ is the  row Brunovsky partition of $H(s)+P(s)$.
We have $\rank(H(s)+P(s))=\rank(G(s))=\rank (H_1(s))+1=\rank( H(s))+1$,
$\widehat s_1-s_1=-1 $ and 
$\widehat s_{100}-s_{100}=99-0=99=-1+\left [ \sqrt{10000}\right ]=-1+\left [ \sqrt{\abs{\bs_*}}\right ]$; i.e. the bounds in (\ref{eqgpboundsr+r}) are reachable.
  
\end{itemize}
\end{enumerate}

\section*{Disclosure statement}

The authors report there are no competing interests to declare.

\section*{Acknowledgments}
The main problem studied in this paper was posed by Prof. Carsten Trunk in the context of linear relations. We would like to thank him for his idea. It has allowed us to state and give a solution to  the problem for matrix pencils.

  \bibliographystyle{acm} 
  \bibliography{references}

  \end{document}